\documentclass{article}
\usepackage{graphicx,amsmath,amsthm,amsopn,amstext,amscd,amsfonts,amssymb}
\usepackage{natbib}

\def\diag{\mathop{\rm diag}\nolimits}
\def\tr{\mathop{\rm tr}\nolimits}
\def\rank{\mathop{\rm rank}\nolimits}
\def\build#1#2#3{\mathrel{\mathop{#1}\limits^{#2}_{#3}}}

\def\vecw{\mathop{\rm vecw}\nolimits}

\def\etr{\mathop{\rm etr}\nolimits}

\renewenvironment{abstract}
                 {\vspace{6pt}
                  \begin{center}
                  \begin{minipage}{5in}
                  \centerline{\textbf{Abstract}}
                  \noindent\ignorespaces
                 }
                 {\end{minipage}\end{center}}
\newtheorem{thm}{\textbf{Theorem}}[section]
\newtheorem{cor}{\textbf{Corollary}}[section]

\theoremstyle{definition}

\newtheorem{rem}{\textbf{Remark}}[section]

\title{\Large \textbf{Generalised shape theory via pseudo-Wishart distribution}}
\author{
  \textbf{Jos\'e A. D\'{\i}az-Garc\'{\i}a} \thanks{Corresponding author\newline
    }\\
  {\normalsize Department of Statistics and Computation} \\
  {\normalsize Universidad Aut\'onoma Agraria Antonio Narro}\\
  {\normalsize 25350 Buenavista, Saltillo, Coahuila, M\'exico} \\
  {\normalsize E-mail: jadiaz@uaaan.mx} \\[2ex]
  \textbf{Francisco J. Caro-Lopera} \\
  {\normalsize Department of Basic Sciences} \\
  {\normalsize Universidad de Medell\'{\i}n} \\
  {\normalsize Carrera 87 No.30-65, of. 5-103}\\
  {\normalsize Medell\'{\i}n, Colombia}\\
  {\normalsize E-mail: fjcaro@udem.edu.co}\\
}
\date{}
\begin{document}
\maketitle

\begin{abstract}
The non isotropic noncentral elliptical shape distributions via pseudo-Wishart
distribution are founded. This way, the classical shape theory is extended to non
isotropic case and the normality assumption is replaced by assuming a elliptical
distribution. In several cases, the new shape distributions are easily computable and
then the inference procedure can be studied under exact densities. An application in
Biology is studied under the classical gaussian approach and two non gaussian models.
\end{abstract}

\section{Introduction}\label{Sec:Introd}

With the introduction of several innovative statistical and mathematical tools for
high-dimensional data analysis, now the classical multivariate analysis have a new
and modern image. Developments as generalised multivariate analysis, latent variable
analysis, DNA microarray data, pattern recognition, multivariate analysis nonlinear,
data mining, manifold learning, shape theory, etc., open a range of potential
applications in many areas of the knowledge.

As consequence of these new statistical and mathematical tools a new theory can be
considere from the conjunction between generalised multivariate analysis and the
statistical shape theory is termed \emph{Generalised Shape Theory}, in which the
methodology developed in the shape theory under Gaussian models is extended to a
general class of distributions, the elliptically contoured densities.

Having this goal in mind, recall that $\mathbf{X}:N\times K$ has a matrix
multivariate elliptically contoured distribution if its  density with respect to the
Lebesgue measure on $\Re^{NK}$ is given by:
$$
  f_{\mathbf{X}}(\mathbf{X})=\frac{1}{|\mathbf{\Sigma}|^{K/2}|\mathbf{\Theta}|^{N/2}}
  h\left\{\tr\left[(\mathbf{X}-\boldsymbol{\mu})'\mathbf{\Sigma}^{-1}(\mathbf{X}-
  \boldsymbol{\mu})\mathbf{\Theta}^{-1}\right]\right\},
$$
where  $\boldsymbol{\mu} :N\times K$, $ \mathbf{\Sigma} :N\times N$,  $
\mathbf{\Theta} :K\times K$, $\mathbf{\Sigma}$ positive definite
($\mathbf{\Sigma}>\mathbf{0}$), $ \mathbf{\Theta}> \mathbf{0}$. The function $h: \Re
\rightarrow [0,\infty)$ is termed the generator function, and it is such that
$\int_{0}^\infty u^{NK-1}h(u^2)du < \infty$. Such a distribution is denoted by
$\mathbf{X}\sim \mathcal{E}_{N\times K}(\boldsymbol{\mu},\mathbf{\Sigma},
\mathbf{\Theta}, h)$, see \citet{fz:90} and \citet{gv:93}. Observe that this class of
matrix multivariate distributions includes Gaussian, Contaminated Normal, Pearson
type II and VI, Kotz, Jensen-Logistic, Power Exponential, Bessel, among other
distributions; whose distributions have tails that are weighted more or less, and/or
distributions with greater or smaller degree of kurtosis than the Gaussian model.

Now, in shape theory context, it is known that the shape of an object is all
geometrical information that remains after filtering out translation, rotation and
scale information of an original figure (represented by a matrix $\mathbf{X}$)
comprised in $N$ landmarks in $K$ dimensions. Hence, we say that two figures,
$\mathbf{X}_{1}:N\times K$ and $\mathbf{X}_{2}:N\times K$ have the same shape if they
are related with a special similarity transformation $\mathbf{X}_{2}=\beta
\mathbf{X}_{1}\mathbf{H}+\mathbf{1}_{N}\boldsymbol{\gamma}'$, where
$\mathbf{H}:K\times K\in SO(K)=\{\mathbf{H} \in \Re^{K \times K}|\mathbf{HH}' =
\mathbf{H}'\mathbf{H} = \mathbf{I}_{K} \mbox{ and } |\mathbf{H}| = + 1\}$ (the
rotation), $\boldsymbol{\gamma}:K\times 1$ (the translation), $\mathbf{1}_{N}:N\times
1,$ $\mathbf{1}_{N}=(1,1,\ldots,1)'$, and $\beta>0$ (the scale). Thus, in this
context, the shape of a matrix $\mathbf{X}$ is all the geometrical information about
$\mathbf{X}$ that is invariant under Euclidean similarity transformations, see
\citet{GM93} and \citet{DM98}.

In the classical statistical shape theory is assumed that $\mathbf{X}$ has the
isotropic matrix multivariate Gaussian distribution with mean
$\boldsymbol{\mu}_{\mathbf{X}}$, see \citet{GM93}, i.e.
$$
  \mathbf{X} \sim \mathcal{N}_{N \times K}(\boldsymbol{\mu}_{\mathbf{X}},
    \sigma^{2}\mathbf{I}_{N} \otimes \mathbf{I}_{K}).
$$
In the context of the \emph{generalised shape theory}, it is assumed that
$$
  \mathbf{X} \sim \mathcal{E}_{N \times K} (\boldsymbol{\mu}_{{}_{\mathbf{X}}},
  \mathbf{\Sigma}_{{}_{\mathbf{X}}} \otimes \mathbf{\Theta},h).
$$
Thus,  two fundamental extensions of classical shape theory are provided, namely:
\begin{itemize}
  \item The generalised theory assumes a matrix multivariate elliptical
    distribution for the landmark data instead of considering a matrix multivariate
    Gaussian distribution.
  \item Also, the usual isotropic Gaussian condition is replaced by assuming a non
  isotropic elliptical model.  Two important advantages are obtained: first, the
  errors are correlated among landmarks, this is considered with the introduction of
  $\mathbf{\Sigma}$, a $N \times N$ definite positive matrix; and second,  the errors
  are correlated among  coordinates of landmarks, this condition is noticed with
  the introduction of $\mathbf{\Theta}$, a $K \times K$ definite positive matrix.
\end{itemize}

The shape coordinates denoted as $\mathbf{u}$ of $\mathbf{X}$ can be constructed by
several ways in terms of QR decomposition, see \citet{GM93}; and singular value
decomposition (SVD), see \citet{g:91}, \citet{lk:93} and \citet{GM93}. For example,
in terms of the QR decomposition, shape coordinates $\mathbf{u}$ of $\mathbf{X}$ are
constructed in several steps summarised in the expression
\begin{equation}\label{eq:QRSteps}
    \mathbf{L}\mathbf{X}\mathbf{\Theta}^{-1/2}=\mathbf{L}\mathbf{Z}=\mathbf{Y}=\mathbf{TH}=
    r\mathbf{WH}=r\mathbf{W}(\mathbf{u})\mathbf{H},
\end{equation}
Observe that $\boldsymbol{\mu}_{{}_{\mathbf{Z}}} = \boldsymbol{\mu}_{{}_{\mathbf{X}}}
\mathbf{\Theta}^{-1/2}$ and the QR shape coordinates of
$\boldsymbol{\mu}_{{}_{\mathbf{Z}}}$ are defined analogously. The matrix
$\mathbf{L}:(N-1)\times N$ has orthonormal rows to $\textbf{1}=(1,\ldots,1)'$.
$\mathbf{L}$ can be a submatrix of the Helmert matrix, for example. Now, let be
$n=\min(N-1,K)$ and $p=\rank \mu$. In (\ref{eq:QRSteps}), $\mathbf{Y} = \mathbf{TH}$
is the QR decomposition, where $\mathbf{T}:(N-1)\times n$ is lower triangular with
$t_{ii}>0$, $i=1,\ldots, \min(n,K-1)$, and $\mathbf{H}:n\times K$, $\mathbf{H}\in
\mathcal{V}_{n,K} = \{\mathbf{H} \in \Re^{n \times K}| \mathbf{HH}' = \mathbf{I}_{n}
\}$, the Stiefel manifold. Note that $\mathbf{T}$ is invariant to translations and
rotations of $\mathbf{Z}$. The matrix $\mathbf{T}$ is referred as the \textit{QR
size-and-shape} and their elements are the QR size-and-shape coordinates of the
original landmark data $\mathbf{Z}$. Typically in shape analysis there are more
landmarks than dimensions ($N>K$). $\mathbf{H}$ acts on the right to transform
$\Re^{K}$ instead of acting on the left as in the multivariate analysis. In our case
we see the landmarks as variables and the dimensions as observations, then the
transposes of our matrices $\mathbf{Z}$ and $\mathbf{Y}$ can be seen as classical
multivariate data matrices. Now, if we divide $\mathbf{T}$ by its size, the centroid
size of $\mathbf{Z}$,
$$
r=\|\mathbf{T}\|=\sqrt{\tr \mathbf{T}'\mathbf{T}}=\|\mathbf{Y}\|.
$$
we obtain the so-termed \textit{QR shape} matrix $\mathbf{W}$ in (\ref{eq:QRSteps}).
Note that $\|\mathbf{W}\|=1$, the elements of $\mathbf{W}$ are a direction vector for
shape, and $\mathbf{u}$ comprises $m=(N-1)K-nK+n(n+1)/2-1$ generalised polar
coordinates.

Observe that, if $\mathbf{\Theta}^{1/2}$ is the positive definite square root of the
matrix $\mathbf{\Theta}$, i .e. $\mathbf{\Theta} = (\mathbf{\Theta}^{1/2})^{2}$, with
$\mathbf{\Theta}^{1/2}:$ $K \times K$, \citet[p. 11]{gv:93}, and noting that
$$
  \mathbf{X} \mathbf{\Theta}^{-1} \mathbf{X}' = \mathbf{X} (\mathbf{\Theta}^{1/2}
  \mathbf{\Theta}^{1/2})^{-1}\mathbf{X}' = \mathbf{X} \mathbf{\Theta}^{-1/2}
  (\mathbf{X} \mathbf{\Theta}^{-1/2})' = \mathbf{Z}\mathbf{Z}',
$$ where
$$
\mathbf{Z} = \mathbf{X} \mathbf{\Theta}^{-1/2},
$$
then
$$
  \mathbf{Z} \sim \mathcal{E}_{N \times K}(\boldsymbol{\mu}_{{}_{\mathbf{Z}}},
  \mathbf{\Sigma}_{{}_{\mathbf{X}}}, \mathbf{I}_{K},h)
$$
with $\boldsymbol{\mu}_{{}_{\mathbf{Z}}} = \boldsymbol{\mu}_{{}_{\mathbf{X}}}
\mathbf{\Theta}^{-1/2}$, see \citet[p. 20]{gv:93}.

And we arrive at the classical starting point in shape theory where the original
landmark matrix is replaced by $\mathbf{Z} = \mathbf{X} \mathbf{\Theta}^{-1/2}$. Then
we can proceed as usual, removing from $\mathbf{Z}$, translation, scale, rotation
and/or reflection in order to obtain the shape of $\mathbf{Z}$ (or $\mathbf{X}$) via
the QR decomposition, for example.

Let be $\boldsymbol{\mu}=\mathbf{L}\boldsymbol{\mu}_{{}_{\mathbf{X}}}$, then
$\mathbf{Y}:(N-1)\times K$ is invariant to translations of the figure $\mathbf{Z}$,
and
$$
  \mathbf{Y}\sim\mathcal{E}_{N-1 \times K}(\boldsymbol{\mu}\mathbf{\Theta}^{-1/2},
  \mathbf{\Sigma}\otimes \mathbf{I}_{K},h),
$$
where $\mathbf{\Sigma}=\mathbf{L}\mathbf{\Sigma}_{\mathbf{X}}\mathbf{L}'$.

As suggest \citet{GM93}, the density of $\mathbf{YY}'$ essentially is the refection
size-and-shape distribution of $\mathbf{Y}$, moreover, it is invariant to orientation
and reflection. Recall that for a given $\mathbf{Y}:N-1\times K$, $n=N-1 < K$, then
$\mathbf{V} = \mathbf{YY}'$ has the noncentral Wishart distribution with respect to
Lebesgue measure on the subspace of definite positive matrices $\mathbf{V} >
\mathbf{0}$. However, the density of $\mathbf{V} = \mathbf{YY}'$ when, $n\geq K$,
exist on the $(nK-K(K-1)/2)$-dimensional manifold of rank-K positive semidefinite
$N-1\times N-1$ matrices with $K$ distinct positive eigenvalues, which is termed
\textit{pseudo-Wishart distribution}, see \citet{u:94}, \citet{dggg05} and
\citet{dggj06}. Therefore, alternatively to (\ref{eq:QRSteps}) we propose the
following steeps for obtain the shape coordinates
\begin{equation}\label{eq:PWSteps}
    \mathbf{L}\mathbf{X}\mathbf{\Theta}^{-1/2}=\mathbf{L}\mathbf{Z}=\mathbf{Y}\Rightarrow\mathbf{V}=
    r\mathbf{W}=r\mathbf{W}(\mathbf{u}),
\end{equation}
where $\mathbf{V} = \mathbf{YY}'$ and $\mathbf{W} = \mathbf{V}/r$, with $r =
||\mathbf{V}||$.

In this work the size and shape distribution for any elliptical model in terms of
pseudo-Wishart distribution is derived in section \ref{sec:PWsizeandshape}. Then the
shape density is obtained in section \ref{sec:PWshape}. The central case of the shape
density is studied in section \ref{sec:PWcentralcase}, and is established that the
central QR reflection shape density is invariant under the elliptical family. Some
particular shape densities are derived in section \ref{sec:particularmodels} in order
to perform inference on exact distributions; i.e. a subfamily of shape distributions
generated by Kotz distributions including the Gaussian is obtained and applied.
Finally in section \ref{sub:mouse}, two elements of that class (the Gaussian and a
non Gaussian model) are applied to an existing publish data, the mouse vertebra
study. Some test for detecting shape differences are gotten and the models are
discriminated by the use of a dimension criterion such as the modified $BIC^{*}$
criterion.

\section{Pseudo-Wishart size-and-shape distribution}\label{sec:PWsizeandshape}

Let $\mathbf{V}= \mathbf{YY}'$. In general ($n=N-1 < K$ or $n\geq K$), the matrix
$\mathbf{V}$ can be written as
$$
   \mathbf{V} \equiv \left (
          \begin{array}{cc}
            \build{\mathbf{V}_{11}}{}{n \times n} & \build{\mathbf{V}_{12}}{}{n \times (N-1)-n} \\
            \build{\mathbf{V}_{21}}{}{(N-1)-n \times n} & \build{\mathbf{V}_{22}}{}{(N-1)-n \times (N-1)-n} \\
          \end{array}
       \right ) \qquad \mbox{with rank of }\mathbf{V}_{11} = n,
$$
such that, the number of mathematically independent elements in
$\mathbf{V}$ are $m=(N-1)K-nK+n(n+1)/2$ corresponding to the mathematically
independent elements in $\mathbf{V}_{11} > \mathbf{0}$ if $n=N-1 < K$ or to the
mathematically independent elements of $\mathbf{V}_{12}$, and
$\mathbf{V}_{11}>\mathbf{0}$ if $n\geq K$. Recall that $\mathbf{V}_{11} >
\mathbf{0}$, in such a way that $\mathbf{V}_{11}$ has $n(n + 1)/2$ mathematically
independent elements, therefore,
\begin{equation*}
    (d\mathbf{V}) \equiv
    \left\{%
    \begin{array}{ll}
        (d\mathbf{V}_{11})=  \displaystyle\bigwedge_{i\leq j}^{n}
        dv_{ij}, & \hbox{if $n=N-1 < K$;} \\
        (d\mathbf{V}_{11})\wedge(d\mathbf{V}_{12})
        = \displaystyle\bigwedge_{i = 1}^{n} \bigwedge_{j = i}^{(N-1)}dv_{ij}, & \hbox{if $n\geq K$.} \\
    \end{array}%
    \right.
\end{equation*}
Formally, the measure $(d\mathbf{V})$ is the Hausdorff measure defined on subspace of
positive semidefinite matrices, see \citet{bi:86}, \citet{dg:97}, \citet{dgm:97} and
\citet{dggg05}.

Explicit forms for $(d\mathbf{V})$ can be obtained under diverse factorisations of
the measure $(d\mathbf{V})$. For example, by using the Cholesky decomposition
$\mathbf{V} = \mathbf{TT}'$, where $\mathbf{T}:(N-1)\times n$ is lower triangular
with $t_{ii}>0$, $i=1,\ldots, \min(n,K-1)$
\begin{equation}\label{mea1}
    (d\mathbf{V})= 2^{n} \displaystyle\prod_{i=1}^{n}t_{ii}^{N-i}(d\mathbf{T}).
\end{equation}
Alternatively, under the nonsingular part of the spectral decompositions $\mathbf{V}
= \mathbf{W}'_{1}\mathbf{D}\mathbf{W}_{1}$, $\mathbf{W}_{1} \in \mathcal{V}_{n,N-1}$
and $\mathbf{D} = \diag(d_{1}, \dots, d_{n})$, $d_{1}> \cdots > d_{n} > 0$ , then
\begin{equation}\label{mea2}
    (d\mathbf{V})= 2^{-n}|\mathbf{D}|^{N-1-n}\displaystyle\prod_{i<j}^{n}(d_{i} - d_{j})
        (d\mathbf{D})(\mathbf{W}'_{1}d\mathbf{W}_{1}).
\end{equation}
Alternative explicit form for $(d\mathbf{V})$ are given in \citet{dggg05}.

\begin{thm}\label{th:PWsizeandshape}
The pseudo-Wishart size-and-shape density is
\begin{equation}\label{eq:PWsizeandshape}
dF_{\mathbf{V}}(\mathbf{V})=\frac{\pi^{nK/2}|\mathbf{V}^{*}|^{(K-N)/2}}
{\Gamma_{n}\left[K/2\right]|\mathbf{\Sigma}|^{K/2}}
\sum_{t=0}^{\infty}\sum_{\kappa}\frac{h^{(2t)}[\tr
(\mathbf{\Sigma}^{-1}\mathbf{V}+\mathbf{\Omega})]}{\left(\frac{1}{2}K\right)_{\kappa}}
\frac{C_{\kappa}(\mathbf{\Omega} \mathbf{\Sigma}^{-1}\mathbf{V})}{t!}(d\mathbf{V}),
\end{equation}
\end{thm}
where $(d\mathbf{V})$ is defined in (\ref{mea1}) or (\ref{mea2}) (among many others),
$\mathbf{\Omega} = \mathbf{\Sigma}^{-1} \boldsymbol{\mu} \mathbf{\Theta}^{-1}
\boldsymbol{\mu}'$, $C_{\kappa}(\mathbf{B})$ are the zonal polynomials of
$\mathbf{B}$ corresponding to the partition $\kappa=(t_{1},\ldots t_{\alpha})$ of
$t$, with $\sum_{i=1}^{\alpha}t_{i}=t$; and
$(a)_{\kappa}=\prod_{i=1}(a-(j-1)/2)_{t_{j}}$, $(a)_{t}=a(a+1)\cdots (a+t-1)$, are
the generalized hypergeometric coefficients and
$\Gamma_{s}(a)=\pi^{s(s-1)/4}\prod_{j=1}^{s}\Gamma(a-(j-1)/2)$ is the multivariate
Gamma function, see \citet{JAT64} and \citet{MR1982}. And $h^{(j)}(v)$ is the $j$-th
derivative of $h$ with respect to $v$. The matrix $\mathbf{V}^{*}$ is given as,
$$
  \mathbf{V}^{*} =
  \left\{%
\begin{array}{ll}
    \mathbf{V}_{11}, & \hbox{under Cholesky decomposition;} \\
    \mathbf{D}, & \hbox{under spectral decomposition.} \\
\end{array}%
\right.
$$
\textit{Proof.} See \citet{dggg05}.

Observe that the density functions (\ref{eq:PWsizeandshape}) with respect to
corresponding Hausdorff measure (\ref{mea1}) or (\ref{mea2}) are not unique,
moreover, the Hausdorff measures (\ref{mea1}) or (\ref{mea2}) are also not unique;
however, from a practical point of view, for example, the maximum likelihood
estimation of the unknown parameters is invariant under different choices of measures
(\ref{mea1}) or (\ref{mea2}) and their corresponding density functions
(\ref{eq:PWsizeandshape}), see \citet[p. 275]{k:68} and \citet[p. 532]{r:73}.

\section{Pseudo-Wishart shape distribution}\label{sec:PWshape}
Observe that for $\mathbf{V}:N-1 \times N-1$, of rank $n=\min (N-1,K)$, hence the
matrix $\mathbf{V}$ contains $(N-1)K-nK+n(n+1)/2$ mathematically independent
pseudo-Wishart coordinates $(v_{ij})$. Let $\vecw \mathbf{V}$ a vector consisting of
mathematically independent elements of $\mathbf{V}$, taken column by column. Then the
pseudo-Wishart shape matrix $\mathbf{W}$ can be written as
$$
\vecw \mathbf{W}=\frac{1}{r}\vecw \mathbf{V},\quad r=||\mathbf{V}||=\sqrt{\tr
\mathbf{V}^{2}} = \sqrt{\tr (\mathbf{Y}'\mathbf{Y})^{2}},
$$
then by  \citet[Theorem 2.1.3, p.55]{MR1982},
$$
  (d \vecw \mathbf{V}) = r^{m}\prod_{i=1}^{m}\sin^{m-i}\theta_{i}
  \left(\bigwedge_{i=1}^{m}d\theta_{i}\right)\wedge dr,
$$
with $m=(N-1)K-nK+n(n+1)/2-1$.  Denoting $\mathbf{u} = (\theta_{1}, \ldots,
\theta_{m})'$, $(d\mathbf{u}) = \bigwedge_{i=1}^{m}d\theta_{i}$ and $J(\mathbf{u}) =
\prod_{i=1}^{m}\sin^{m-i}\theta_{i}$, with  $r>0$, $0 < \theta_{i} \leq \pi$ ($i = 1,
\dots, m-1$), $0 < \theta_{m} \leq 2\pi$, then
$$
(d\mathbf{V}) = r^{m} J(\mathbf{u})(d\mathbf{u})\wedge dr.
$$

\begin{thm}\label{th:PWreflectionshape}
The pseudo-Wishart reflection shape density is
\begin{eqnarray}\label{eq:QRreflectionshape}
dF_{\mathbf{W}}(\mathbf{W})&=&\frac{\pi^{nK/2}|\mathbf{W}^{*}|^{(K-N)/2}J(\mathbf{u})}
{\Gamma_{n}\left[K/2\right]|\mathbf{\Sigma}|^{K/2}}
\sum_{t=0}^{\infty}\sum_{\kappa}\frac{C_{\kappa}(\mathbf{\Omega}
\mathbf{\Sigma}^{-1}\mathbf{W})}{t!\left(\frac{1}{2}K\right)_{\kappa}}\nonumber\\
&&\times\int_{0}^{\infty} r^{m-n(K-N)/2+t}h^{(2t)}[r\tr
\mathbf{\Sigma}^{-1}\mathbf{W}+\tr\mathbf{\Omega}](dr)(d\mathbf{u}),
\end{eqnarray}
\end{thm}
where $\mathbf{W}^{*} = \mathbf{V}^{*}/r$.

\textit{Proof.}  The density of $\mathbf{V}$ is
\begin{equation*}
dF_{\mathbf{V}}(\mathbf{V})=\frac{\pi^{nK/2}|\mathbf{V}^{*}|^{(K-N)/2}}
{\Gamma_{n}\left[K/2\right]|\mathbf{\Sigma}|^{K/2}}
\sum_{t=0}^{\infty}\sum_{\kappa}\frac{h^{(2t)}[\tr
(\mathbf{\Sigma}^{-1}\mathbf{V}+\mathbf{\Omega})]}{\left(\frac{1}{2}K\right)_{\kappa}}
\frac{C_{\kappa}(\mathbf{\Omega} \mathbf{\Sigma}^{-1}\mathbf{V})}{t!}(d\mathbf{V}).
\end{equation*}
Making the change of variables $\mathbf{W}(\mathbf{u})=\mathbf{V}/r$, the joint
density function of $r$ and $\mathbf{u}$ is
\begin{eqnarray*}
f_{r,\mathbf{W}}(r,\mathbf{W})&=&\frac{\pi^{nK/2}|r\mathbf{W}^{*}|^{(K-N)/2}}
{\Gamma_{n}\left[K/2\right]|\mathbf{\Sigma}|^{K/2}}
\sum_{t=0}^{\infty}\sum_{\kappa}\frac{h^{(2t)}[\tr (r\mathbf{\Sigma}^{-1}\mathbf{W}
+\mathbf{\Omega})]}{\left(\frac{1}{2}K\right)_{\kappa}}\\
&&\times\frac{C_{\kappa}(r\mathbf{\Omega}\mathbf{\Sigma}^{-1}\mathbf{W})}
{t!}r^{m}J(\mathbf{u})dr \wedge (d\mathbf{u}).
\end{eqnarray*}
Now, note that
\begin{itemize}
    \item
    $C_{\kappa}(r\mathbf{\Omega}\mathbf{\Sigma}^{-1}\mathbf{W})=r^{t}
    C_{\kappa}(\mathbf{\Omega}\mathbf{\Sigma}^{-1}\mathbf{W})$.
    \item $
    |r\mathbf{W}^{*}|^{(K-N)/2} = r^{n(K-N)/2} |\mathbf{W}^{*}|^{(K-N)/2}.
    $
    \item
    $h^{(2t)}[\tr(r\mathbf{\Sigma}^{-1}\mathbf{W}+\mathbf{\Omega})] =
    h^{(2t)}[r\tr\mathbf{\Sigma}^{-1}\mathbf{W}+\tr\mathbf{\Omega}]$.
\end{itemize}
Finally, collecting powers of $r$ by $r^{m+n(K-N)/2 + t}$,  the marginal of
$\mathbf{W}$ is obtained integrating with respect to $r$. \qed

When $\mathbf{\Sigma} = \sigma^{2}\mathbf{I}$, then $\mathbf{\Omega} =
\boldsymbol{\mu} \mathbf{\Theta}^{-1} \boldsymbol{\mu}'/\sigma^{2}$,
$|\mathbf{\Sigma}|^{K/2} = \sigma^{M}$, $M = (N-1)K$ and $r\tr
\mathbf{\Sigma}^{-1}\mathbf{W} = r\tr \mathbf{W}/\sigma^{2}$, thus Theorem
\ref{th:PWreflectionshape} becomes.

\begin{cor}\label{coroSHdiagonal}
The isotropic pseudo-Wishart reflection shape density is
\begin{eqnarray}\label{eq:PWreflectionshapediagonal}
dF_{\mathbf{W}}(\mathbf{W})&=&\frac{\pi^{nK/2}|\mathbf{W}^{*}|^{(K-N)/2}J(\mathbf{u})}
{\Gamma_{n}\left[\frac{1}{2}K\right]\sigma^{M}}
\sum_{t=0}^{\infty}\sum_{\kappa}\frac{C_{\kappa}\left(\displaystyle\frac{1}{\sigma^{2}}\mathbf{\Omega}
\mathbf{W}\right)}{t!\left(\frac{1}{2}K\right)_{\kappa}}\nonumber\\
&& \times \int_{0}^{\infty}
r^{m-n(K-N)/2+t}h^{(2t)}[r\tr\mathbf{W}/\sigma^{2}+\tr\mathbf{\Omega}](dr)\wedge
(d\mathbf{u}).
\end{eqnarray}
\end{cor}

\section{Central case}\label{sec:PWcentralcase}
The central case of the preceding sections can be derived easily.
\begin{cor}\label{cor:PWcentralreflectionsizeandshape}
The central pseudo-Wishart reflection size-and-shape density is given by
\begin{equation*}
dF_{\mathbf{V}}(\mathbf{V})=\frac{\pi^{nK/2}|\mathbf{V}^{*}|^{(K-N)/2}}
{\Gamma_{n}\left[K/2\right]|\mathbf{\Sigma}|^{K/2}} h[\tr
\mathbf{\Sigma}^{-1}\mathbf{V}](d\mathbf{V}).
\end{equation*}
\end{cor}
\textit{Proof.} It is straightforward from Theorem \ref{th:PWsizeandshape} just take
$\boldsymbol{\mu} = \mathbf{0}$ and recall that $h^{(0)}[\tr\cdot]=h[\tr\cdot]$. \qed

Similarly:
\begin{cor}\label{cor:PWcentralreflectionshape}
The central pseudo-Wishart reflection shape density is given by
$$
dF_{\mathbf{W}}(\mathbf{W})=\frac{\pi^{nK/2}|\mathbf{W}^{*}|^{(K-N)/2}J(\mathbf{u})}
{\Gamma_{n}\left[K/2\right]|\mathbf{\Sigma}|^{K/2}} \int_{0}^{\infty}
r^{m-n(K-N)/2}h[r\tr \mathbf{\Sigma}^{-1}\mathbf{W}](dr)(d\mathbf{u}).
$$
\end{cor}
\textit{Proof.} Just take $\boldsymbol{\mu} = \mathbf{0}$ and
$h^{(0)}[\tr\cdot]=h[\tr\cdot]$ in Theorem \ref{th:PWreflectionshape}. \qed

Observe that it is possible to obtain an invariant central shape density, i.e. the
density  does not depend on function $h(\cdot)$ Let $h$ be the density generator of
$\mathbf{Y}\sim \mathcal{E}_{N-1,K}(\mathbf{0},\mathbf{I} \otimes \mathbf{I},h)$,
i.e.
$$
f_{\mathbf{Y}}(\mathbf{Y})=h(\tr \mathbf{Y}\mathbf{Y}'),
$$
then by \citet[eq. 3.2.6, p.102]{fz:90},
$$
  \int_{0}^{\infty}r^{(N-1)K-1}h(r^{2})dr=\frac{\Gamma[(N-1)K/2]}{2\pi^{(N-1)K/2}}.
$$
Taking $s = r^{2}$ with $dr = ds/(2\sqrt{s})$
$$
  \int_{0}^{\infty}s^{(N-1)K/2-1}h(s)dr=\frac{\Gamma[(N-1)K/2]}{\pi^{(N-1)K/2}}.
$$
Hence, if $s=(\tr\mathbf{\Sigma}^{-1}\mathbf{W})r$,
$ds=(\tr\mathbf{\Sigma}^{-1}\mathbf{W})(dr)$, then

$
  \displaystyle\int_{0}^{\infty}r^{m-n(K-N)/2}h[r\tr\mathbf{\Sigma}^{-1}\mathbf{W}](dr)
$
\begin{eqnarray*}
&=&\int_{0}^{\infty}\left(\frac{s}{(\tr\mathbf{\Sigma}^{-1}\mathbf{W})}\right)^{m-n(K-N)/2}
h(s)\frac{ds}{(\tr\mathbf{\Sigma}^{-1}\mathbf{W})}\\
&=&(\tr\mathbf{\Sigma}^{-1}\mathbf{W})^{n(K-N)/2-m-1}\int_{0}^{\infty}s^{(2m-n(K-N))/2+1-1}h(s)ds\\
&=&(\tr\mathbf{\Sigma}^{-1}\mathbf{W})^{n(K-N)/2-m-1}\frac{\Gamma[m-n(K-N)/2+1]}{\pi^{m-n(K-N)/2+1}}.
\end{eqnarray*}
Thus:
\begin{cor}\label{cor:PWcentralreflectionshapeinvariant}
When $\boldsymbol{\mu}=\mathbf{0}$ the pseudo-Wishart reflection shape density is
invariant under the elliptical family and it is given by
\begin{eqnarray*}
    dF_{\mathbf{W}}(\mathbf{W})&=&\frac{\pi^{nK-m+n(K-N)/2-1}\Gamma[m-n(K-N)/2+1]}{\Gamma_{n}
    \left[K/2\right]|\mathbf{\Sigma}|^{K/2}}
   |\mathbf{W}^{*}|^{(K-N)/2}\\
   && \quad \times J(\mathbf{u})(\tr\mathbf{\Sigma}^{-1}\mathbf{W})^{n(K-N)/2-m-1}(d\mathbf{u}).
\end{eqnarray*}
\end{cor}
Now, if $\mathbf{\Sigma} = \sigma^{2}\mathbf{I}$, then
$$
  (\tr\mathbf{\Sigma}^{-1}\mathbf{W})^{n(K-N)/2-m-1} = (1/\sigma^{2})^{n(K-N)/2-m-1}
  (\tr \mathbf{W})^{n(K-N)/2-m-1},
$$
and $|\mathbf{\Sigma}|^{K/2} = (\sigma^{2})^{M/2}$, thus:
\begin{cor}\label{cor:QPWcentralreflectionshapeinvariantdiagonal}
When $\boldsymbol{\mu}=\mathbf{0}$ and $\mathbf{\Sigma} = \sigma^{2}\mathbf{I}$ the
pseudo-Wishart reflection shape density is invariant under the elliptical family and
it is given by
\begin{eqnarray*}
   dF_{\mathbf{W}}(\mathbf{W})&=&\frac{\pi^{nK-m+n(K-N)/2-1}\Gamma[m-n(K-N)/2+1]}{2\Gamma_{n}
    \left[K/2\right](\sigma^{2})^{n(K-N)/2+M/2-m-1}}
   |\mathbf{W}^{*}|^{(K-N)/2}\\
   && \quad (\tr \mathbf{W})^{n(K-N)/2-m-1}\times J(\mathbf{u})(d\mathbf{u}).
\end{eqnarray*}
\end{cor}

\section{Some particular models}\label{sec:particularmodels}

Finally, we give explicit shapes densities for some elliptical models.

The Kotz type I model is given by
\begin{equation}\label{Kotz1}
    h(y) = \frac{R^{T-1+\frac{K(N-1)}{2}}\Gamma\left[\frac{K(N-1)}{2}\right]}{\pi^{K(N-1)/2}
  \Gamma\left[T-1+\frac{K(N-1)}{2}\right]}y^{T-1}\exp\{-Ry\},
\end{equation}
Then, the corresponding $k$-th derivative $\displaystyle
\frac{d^{k}[y^{T-1}\exp\{-Ry\}]}{dy^{k}}$, is
\begin{equation}\label{eq:kotzcase1}
\frac{(-R)^{k}y^{T-1}}{\exp\{Ry\}}\left\{1+\sum_{m=1}^{k}\binom{k}{m}
\left[\prod_{i=0}^{m-1}(T-1-i)\right](-Ry)^{-m}\right\}.
\end{equation}

It includes the Gaussian case, i.e. when $T=1$ and $R=1/2$, here the derivation is
straightforward from the general density.

The required derivative follows easily, it is,
$$
  h^{(k)}(y)=\frac{R^{M/2}}{\pi^{M/2}}(-R)^{k}\exp(-Ry)
$$
and
\begin{small}
\begin{eqnarray*}
   \int_{0}^{\infty}&&r^{m-n(K-N)/2+t}h^{(2t)}[r\tr
\mathbf{\Sigma}^{-1}\mathbf{W}+\tr\mathbf{\Omega}]dr\\&&=
\pi^{-M/2}R^{-m+t+\frac{1}{2}(-2+M+n(K-N))}(\tr
\mathbf{\Sigma}^{-1}\mathbf{W})^{-1-m-t+n(K-N)/2}\\&&
\times\etr\left(-R\mathbf{\Omega}\right)\Gamma\left[1+m+t+\frac{1}{2}n(-K+N)\right].
\end{eqnarray*}
\end{small}
So, we have proved that

\begin{cor}\label{cor:reflectionshapeRNORMAL}
The Kotz type I ($T=1$) Pseudo-Wishart reflection shape density  is
\begin{eqnarray*}
dF_{\mathbf{W}}(\mathbf{W})&=&
\frac{\pi^{(nK-M)/2}|\mathbf{W}^{*}|^{(K-N)/2}J(\mathbf{u})\etr\left(-R\mathbf{\Omega}\right)}
{R^{m-\frac{1}{2}(-2+M+n(K-N))}\Gamma_{n}\left[K/2\right]|\mathbf{\Sigma}|^{K/2}}
\nonumber\\
&& \ \times \sum_{t=0}^{\infty}\frac{\Gamma\left[1+m+t+n(-K+N)/2\right]}{t!(\tr
\mathbf{\Sigma}^{-1}\mathbf{W})^{1+m+t-n(K-N)/2}}\sum_{\kappa}\frac{C_{\kappa}(R\mathbf{\Omega}
\mathbf{\Sigma}^{-1}\mathbf{W})}{\left(\frac{1}{2}K\right)_{\kappa}}.
\end{eqnarray*}
where $M=(N-1)K$.
\end{cor}

Finally, for the  Kotz type I model  (\ref{Kotz1}) and the given $2t$-th derivative,
we can prove easily  that

\begin{cor}\label{cor:SVDreflectionshapeKotz}
The pseudo-Wishart reflection shape density based on the Kotz type I model is given
by
\begin{equation}\label{eq:coroQRreflectionshape}
dF_{\mathbf{W}}(\mathbf{W})=\frac{\pi^{nK/2}|\mathbf{W}^{*}|^{(K-N)/2}J(\mathbf{u})}
{\Gamma_{n}\left[K/2\right]|\mathbf{\Sigma}|^{K/2}}
\sum_{t=0}^{\infty}\sum_{\kappa}\frac{C_{\kappa}(\mathbf{\Omega}
\mathbf{\Sigma}^{-1}\mathbf{W})}{t!\left(\frac{1}{2}K\right)_{\kappa}}I(\mathbf{W(u)},r)
\ (d\mathbf{u})
\end{equation}
where
\begin{small}
\begin{eqnarray*}
&&I(\mathbf{W(u)},r) =\int_{0}^{\infty}
r^{m-n(K-N)/2+t}h^{(2t)}[r\tr
\mathbf{\Sigma}^{-1}\mathbf{W}+\tr\mathbf{\Omega}](dr)(d\mathbf{u})\\
    &=&G\
    e^{-RB}A^{-a-1}\left[\sum_{u=0}^{\infty}(u!)^{-1}R^{2t-1-a-u}B^{T-1-u}\Gamma[1+a+u]\prod_{s=0}^{u-1}(T-1-s)\right.
    \\&&\left.+\sum_{v=1}^{2t}\binom{2t}{v}\left[\prod_{i=0}^{v-1}(T-1-i)\right]
    \sum_{u=0}^{\infty}(u!)^{-1}(-1)^{-v}R^{2t-1-a-u-v}B^{T-1-u-v}\right.\\&&\left.\times\Gamma[1+a+u]\prod_{s=0}^{u-1}(T-1-v-s)\right],
\end{eqnarray*}
\end{small}
with $M = (N-1)K$, $G = \pi^{-M/2}R^{T-1+M/2}\Gamma[M/2]/
\Gamma\left[T-1+M/2\right]$, $A = \tr \mathbf{\Sigma}^{-1}\mathbf{W}$, $B =
\tr\boldsymbol{\Omega}$ and $a = m-n(K-N)/2+t$.
\end{cor}
This density seems uncomputable but it easy to see that it has the form of a
generalised hypergeometric functions (see next section). These series can be
determined by suitable modifications of the algorithms given by \citet{KE06} for
${}_{0}F_{1}$ and at the same computational costs. Moreover, if the parameter $T>0$
is an integer,  the series are simplified substantially. For example,  we can prove
that the shape density associated to a Kotz model with $T=3$, $R=1/2$ and the
isotropic assumption ($\mathbf{\Sigma}=\sigma^{2}\mathbf{I}_{N-1}$ and
$\mathbf{\Theta}=\mathbf{I}_{K}$), is given by:
\begin{small}
\begin{eqnarray}\label{eq:pwT3}
&&dF_{\mathbf{W}}(\mathbf{W})=\frac{\pi^{(nK-M)/2}|\mathbf{W}^{*}|^{(K-N)/2}J(\mathbf{u})\etr\left(-\boldsymbol{\mu}'
\boldsymbol{\mu}/2\sigma^{2}\right)}
{2^{-3-m+(M+n(K-N))/2}M(M+2)\Gamma_{n}\left[K/2\right]}\\
&&\times
\sum_{t=0}^{\infty}\frac{[(B-2t)^{2}-2t]\Gamma\left[a\right]+2(B-2t)\Gamma\left[a+1\right]+\Gamma\left[a+2\right]}
{t!\sigma^{M-2-2m+n(K-N)}\left(\tr\mathbf{W}\right)^{a}}\sum_{\kappa}
\frac{C_{\kappa}(\frac{1}{2\sigma^{2}}
\boldsymbol{\mu}'\mathbf{W}\boldsymbol{\mu})}{\left(\frac{1}{2}K\right)_{\kappa}}\nonumber.
\end{eqnarray}
\end{small}
where $M=(N-1)K$, $B=\tr\boldsymbol{\mu}'\boldsymbol{\mu}/2\sigma^{2}$ and
$a=1+m+t+n(-K+N)/2$.

Other examples shall be considered in the next section, when $T=1$ and $T=2$. More
complex densities in the context of affine shape theory were computed by using the
same idea, see \citet{Caro2009}.

\section{Example}\label{sub:mouse}

\begin{figure}[ht]
  \begin{center}
   \includegraphics[width=7cm,height=7cm]{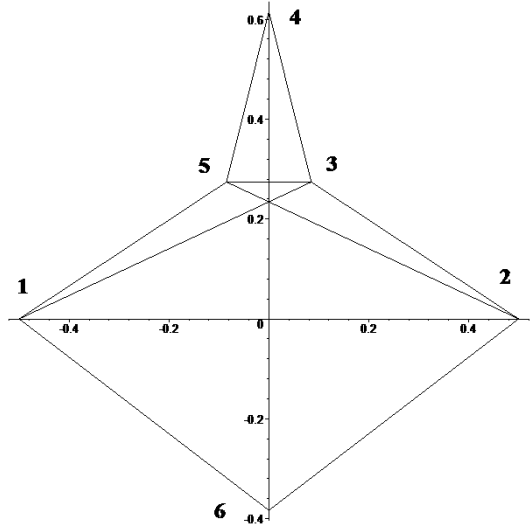}
   \caption{Mouse vertebra}\label{fig:mouse}
  \end{center}
\end{figure}

This problem is studied in detail by \citet{DM98} under a number of approaches (see
also \citet{MD:89}). The experiment considers the second thoracic vertebra T2 of two
groups of mice: large and small. The mice are selected and classified according to
large or small body weight, respectively; in this case, the sample consists of 23
large and small bones (the data can be found in \citet[p. 313-316]{DM98}). It is of
interest to study shape differences between the two groups. The vertebras are
digitised and summarised in six mathematical landmarks which are placed at points of
high curvature, see figure \ref{fig:mouse}; they are symmetrically selected by
measuring the extreme positive and negative curvature of the bone. See \citet{DM98}
for more details.

\medskip

Here we study three models, the  Gaussian shape, and two shape Kotz type I models
with $T=2$ and $T=3$.

First, the isotropic Gaussian shape density is obtained from corollary
\ref{cor:reflectionshapeRNORMAL} when we set $R=\frac{1}{2}$,
$\mathbf{\Sigma}=\sigma^{2}I_{N-1}$, $\mathbf{\Theta}=I_{K}$, $\mathbf{\Omega} =
\mathbf{\Sigma}^{-1} \boldsymbol{\mu} \mathbf{\Theta}^{-1} \boldsymbol{\mu}'=
\sigma^{-2} \boldsymbol{\mu} \boldsymbol{\mu}'$, namely

\begin{cor}\label{cor:isotropicnormal}
The Pseudo-Wishart reflection shape density based on the isotropic Gaussian is given by
\begin{small}
\begin{eqnarray*}
dF_{\mathbf{W}}(\mathbf{W})&=&\frac{\pi^{(nK-M)/2}|\mathbf{W}^{*}|^{(K-N)/2}J(\mathbf{u})\etr\left(-\boldsymbol{\mu}'
\boldsymbol{\mu}/2\sigma^{2}\right)}
{2^{-m+(-2+M+n(K-N))/2}\Gamma_{n}\left[K/2\right]}
\nonumber\\
&&\times
\sum_{t=0}^{\infty}\frac{\Gamma\left[1+m+t+n(-K+N)/2\right]}{t!\sigma^{M-2-2m+n(K-N)}
\left(\tr\mathbf{W}\right)^{ m-n(K-N)/2+t+1}}\sum_{\kappa}
\frac{C_{\kappa}(\frac{1}{2\sigma^{2}}
\boldsymbol{\mu}'\mathbf{W}\boldsymbol{\mu})}{\left(\frac{1}{2}K\right)_{\kappa}}.
\end{eqnarray*}
\end{small}
where $M=(N-1)K$.
\end{cor}

A second shape distribution that we will use follows from corollary
\ref{cor:SVDreflectionshapeKotz} by taking $R=1/2$, $T=2$, i.e.
\begin{eqnarray*}
&&dF_{\mathbf{W}}(\mathbf{W})=\frac{\pi^{(nK-M)/2}|\mathbf{W}^{*}|^{(K-N)/2}J(\mathbf{u})\etr\left(-\boldsymbol{\mu}'
\boldsymbol{\mu}/2\sigma^{2}\right)}
{2^{-2-m+(M+n(K-N))/2}M\Gamma_{n}\left[K/2\right]}
\nonumber\\
&&\times
\sum_{t=0}^{\infty}\frac{(B-2t)\Gamma\left[a\right]+\Gamma\left[a+1\right]}{t!\sigma^{M-2-2m+n(K-N)}
\left(\tr\mathbf{W}\right)^{a}}\sum_{\kappa} \frac{C_{\kappa}(\frac{1}{2\sigma^{2}}
\boldsymbol{\mu}'\mathbf{W}\boldsymbol{\mu})}{\left(\frac{1}{2}K\right)_{\kappa}}.
\end{eqnarray*}
where $M=(N-1)K$, $B=\tr\boldsymbol{\mu}'\boldsymbol{\mu}/2\sigma^{2}$ and
$a=1+m+t+n(-K+N)/2$.

And the third shape model of this example corresponds to the isotropic Kotz
distribution with $T=3$ and $R=1/2$, see (\ref{eq:pwT3}).

In order to select the best elliptical model, a number of dimension criteria have
been proposed. We shall consider a modification of the $BIC^{*}$ statistic as
discussed in \citet{YY07}, and which was first achieved by \citet{ri:78} in a coding
theory framework. The modified $BIC^{*}$ is given by:
$$
  BIC^{*}=-2\mathfrak{L}(\widetilde{\boldsymbol{\mu}},\widetilde{\sigma}^{2},h)
         + n_{p}(\log(n+2 )- \log 24),
$$
where $\mathfrak{L}(\widetilde{\boldsymbol{\mu}},\widetilde{\sigma}^{2},h)$ is the
maximum of the log-likelihood function, $n$ is the sample size and $n_{p}$ is the
number of parameters to be estimated for each particular shape density.

Now, if the goal of the shape analysis searches the best elliptical distribution,
among a set of proposed models, the modified $BIC^{*}$ criterion suggests to choose
the model for which the modified $BIC^{*}$  receives its smallest value. In addition,
as proposed by \citet{kr:95} and \citet{r:95}, the following selection criteria have
been employed in order to compare two contiguous models in terms of its corresponding
modified $BIC^{*}$ .

\medskip

\begin{table*}[ht]
\centering \caption{Grades of evidence corresponding to values of the $BIC^{*}$
difference.}\label{table2}
\renewcommand{\arraystretch}{1}
\begin{center}
  \begin{tabular}{cl}
    \hline
    $BIC^{*}$ difference & Evidence\\
    \hline
    0--2 & Weak\\
    2--6 & Positive \\
    6--10 & Strong\\
    $>$ 10 & Very strong\\
    \hline
  \end{tabular}
\end{center}
\end{table*}

Now, recall that for a general density generator $h(\cdot)$
$$
  \mathbf{X} \sim \mathcal{E}_{N \times K} (\boldsymbol{\mu}_{{}_{\mathbf{X}}},
  \mathbf{\Sigma}_{{}_{\mathbf{X}}} \otimes \mathbf{\Theta},h),
$$
where
$\boldsymbol{\mu}=\mathbf{L}\boldsymbol{\mu}_{{}_{\mathbf{X}}}$,
then
$$
  \mathbf{Y}\sim\mathcal{E}_{N-1 \times K}(\boldsymbol{\mu}\mathbf{\Theta}^{-1/2},
  \mathbf{\Sigma}\otimes \mathbf{I}_{K},h),
$$
with
$\mathbf{\Sigma}=\mathbf{L}\mathbf{\Sigma}_{\mathbf{X}}\mathbf{L}'$.

In the mouse vertebra experiment, we want to find the maximum likelihood estimators
(MLE) of the mean shape
\begin{equation*}
   \boldsymbol{\mu}=\left(%
\begin{array}{cc}
  \mu_{11} & \mu_{12} \\
  \mu_{21} & \mu_{22} \\
  \mu_{31} & \mu_{32} \\
  \mu_{41} & \mu_{42} \\
  \mu_{51} & \mu_{52} \\
\end{array}%
\right),
\end{equation*}
and the scale parameter $\sigma^{2}$ defined in the isotropy assumption
$\boldsymbol{\Theta}=\mathbf{I}_{K}$ and
$\boldsymbol{\Sigma}=\boldsymbol{\sigma}^{2}\mathbf{I}_{N-1}$, (in order to
accelerate the computations of this example we fix the variance of the process as 50
-the maximum median variance of the two samples-). This optimisation is applied in
the two independent populations, the small and large groups; first by assuming a
Gaussian model and afterwards by considering two Kotz models indexed by $T=2$ and
$T=3$.

The general procedure is the following: Let
$\mathfrak{L}(\widetilde{\boldsymbol{\mu}},\widetilde{\sigma}^{2},h)$ be the log
likelihood function of a given group-model. The maximisation of the likelihood
function $\mathfrak{L}(\widetilde{\boldsymbol{\mu}},\widetilde{\sigma}^{2},h)$, is
obtained in this paper by using the \emph{Nelder-Mead Simplex Method}, which is an
unconstrained multivariable function using a derivative-free method; specifically, we
apply  the routine \textbf{fminsearch} implemented by the sofware MatLab.

As the reader can check, the shape densities are series of zonal polynomials of the
form
\begin{equation}\label{suma}
    \sum_{t=0}^{\infty}\frac{f(t,\tr \mathbf{X} )}{t!}\sum_{\kappa}\frac{C_{\kappa}(\mathbf{X})}{(a)_{\kappa}},
\end{equation}
which has hypergeometric series
$$
  \sum_{t=0}^{\infty}\frac{1}{t!}\sum_{\kappa}\frac{C_{\kappa}(\mathbf{X})}{(a)_{\kappa}},
$$
as a particular case; these series were non computable for decades. The work of
\citet{KE06} solved the problem and it let the computation of the hypergeometric
series by truncation of the series until the coefficient for large degrees are zero
under certain tolerance. The cited algorithm gives the coefficients of the series,
then, we can modified the algorithm for hypergeometric series to compute the shape
densities with the same computational costs, multiplying each coefficient of the
series by the required function $f(t,\tr \mathbf{X})$.

\medskip

\begin{table}[!ht]
  \centering
  \caption{The maximum likelihood estimators for the \textit{small} group under the
  \textit{Gaussian} model}\label{tab:Gaussiansmall}
\begin{center}
\begin{small}
\renewcommand{\arraystretch}{1}
\begin{tabular}{c|c|c|c|c|c|c|c}
  \hline
  % after \\: \hline or \cline{col1-col2} \cline{col3-col4} ...
  Trunc. & $\widetilde{\mu}_{11}$  & $\widetilde{\mu}_{12}$ & $\widetilde{\mu}_{21}$ & $\widetilde{\mu}_{22}$
  & $\widetilde{\mu}_{31}$ & $\widetilde{\mu}_{32}$ & $\widetilde{\mu}_{41}$ \\
  \hline
  20  &30.40  & -8.13 &  5.73  &  9.47  &  4.01  &  17.34   &-2.70\\
  40  & -0.47 &  -44.69 & 15.04   &-4.54 &  25.27 &  0.60  &  4.88\\
  60  &-2.10  & -54.84&  18.31&   -6.09&   31.03 &  -0.12&   6.17\\
  80  &-0.70 &  -63.48 & 21.37  & -6.46&   35.89   &0.83  &  6.94\\
  100 &-2.61  & -71.03 & 23.73 &  -7.85   &40.19 &  -0.10&   7.98\\
  110 &-0.54 &  -74.58 & 25.13 &  -7.50   &42.16  & 1.14 &   8.12\\
  120 & -3.41 &  -77.44 & 25.88 &  -8.71 &  43.94 &  -0.37  & 8.79\\
  140 &   -3.41 &  -77.44 & 25.88 &  -8.71 &  43.94 &  -0.37  & 8.79\\
  160 & -3.41 &  -77.44 & 25.88 &  -8.71 &  43.94 &  -0.37  & 8.79\\
 \hline
\end{tabular}

\medskip

\begin{tabular}{c|c|c|c|c|c|c}
  \hline
  % after \\: \hline or \cline{col1-col2} \cline{col3-col4} ...
  Trunc. & $\widetilde{\mu}_{42}$ & $\widetilde{\mu}_{51}$ & $\widetilde{\mu}_{52}$  & $BIC^{*}$ & Time & Iter. \\
  \hline
  20  &4.24 &   -6.82 &  -20.65 & -3538.26 &   317& 4103\\
  40  &5.21 &   -30.81 & 2.11   & -4155.34  &  281 &1881\\
  60  & 6.23&    -37.75&  3.63 &   -4659.16  &  417 &1923\\
  80  &7.40 &   -43.77&  3.01   & -5110.98 &   426 &1455\\
  100 &8.08 &   -48.90&  4.63  &  -5532.79 &   742& 2025\\
  110 &8.72 &   -51.43 & 3.34  &  -5735.95  &  607 &1507\\
  120 &8.76&    -53.43 & 5.37 &   -5914.74&    721& 1640\\
  140 & 8.76&    -53.43 & 5.37 &   -5914.74&    721& 1640\\
  160 & 8.76&    -53.43 & 5.37 &   -5914.74&    721& 1640\\
 \hline
\end{tabular}
\end{small}
\end{center}
\end{table}

\medskip

\begin{table}[!ht]
  \centering
  \caption{The maximum likelihood estimators for the \textit{small} group under the
  \textit{Kotz $T=2$} model}\label{tab:T2small}
\begin{center}
\begin{small}
\renewcommand{\arraystretch}{1}
\begin{tabular}{c|c|c|c|c|c|c|c}
  \hline
  % after \\: \hline or \cline{col1-col2} \cline{col3-col4} ...
  Trunc. & $\widetilde{\mu}_{11}$  & $\widetilde{\mu}_{12}$ & $\widetilde{\mu}_{21}$ & $\widetilde{\mu}_{22}$
  & $\widetilde{\mu}_{31}$ & $\widetilde{\mu}_{32}$ & $\widetilde{\mu}_{41}$ \\
  \hline
  20  &-6.06  & -32.06 & 10.23 &  -5.19 &  18.24  & -2.80 &  4.18\\
  40  & 4.42&    -46.00&  15.97 &  -3.02  & 25.91 &  3.39  &  4.45\\
  60  &-0.53&   -56.62 & 19.07   &-5.73 &  32.01  & 0.80 &   6.18\\
  80  &-1.88 &  -65.36&  21.89  & -7.04 &  36.97  & 0.20    &7.28\\
  100 &-0.04 &  -73.09 & 24.68  & -7.18   &41.31 &  1.40   & 7.90\\
  110 &-1.62 &  -76.63 & 25.72 &  -8.06 &  43.34   &0.57  &  8.47\\
  120 &-1.84 &  -79.60 & 26.76   &-8.39&   45.13   &0.56  &  8.83\\
  140 &-1.84 &  -79.60 & 26.76   &-8.39&   45.13   &0.56  &  8.83\\
  160 &-1.84 &  -79.60 & 26.76   &-8.39&   45.13   &0.56  &  8.83\\
 \hline
\end{tabular}

\medskip

\begin{tabular}{c|c|c|c|c|c|c}
  \hline
  % after \\: \hline or \cline{col1-col2} \cline{col3-col4} ...
  Trunc. & $\widetilde{\mu}_{42}$ & $\widetilde{\mu}_{51}$ & $\widetilde{\mu}_{52}$  & $BIC^{*}$ & Time & Iter. \\
  \hline
  20  & 3.12 &   -21.88 & 5.46 &   -3584.58 &   311& 1957\\
  40  &5.89&    -31.91 & -1.22  & -4203.67 &   627 &2052\\
  60  &6.61 &   -39.04 & 2.62  &  -4709.27 &   890 &1986\\
  80  &7.49   & -45.02&  3.90 &   -5162.67  &  951& 1566\\
  100 &8.60 &   -50.43 & 2.94   & -5585.92 &   1468    &1978\\
  110 & 8.84   & -52.81 & 4.17   & -5789.74&    1160   & 1386\\
  120 &9.18  &  -54.98 & 4.37  &  -5969.11&    1464&    1656\\
  140 &9.18  &  -54.98 & 4.37  &  -5969.11&    1464&    1656\\
  160 &9.18  &  -54.98 & 4.37  &  -5969.11&    1464&    1656\\
 \hline
\end{tabular}
\end{small}
\end{center}
\end{table}

\medskip

\begin{table}[!ht]
  \centering
  \caption{The maximum likelihood estimators for the \textit{small} group under the
  \textit{Kotz $T=3$} model}\label{tab:T3small}
\begin{center}
\begin{small}
\renewcommand{\arraystretch}{1}
\begin{tabular}{c|c|c|c|c|c|c|c}
  \hline
  % after \\: \hline or \cline{col1-col2} \cline{col3-col4} ...
  Trunc. & $\widetilde{\mu}_{11}$  & $\widetilde{\mu}_{12}$ & $\widetilde{\mu}_{21}$ & $\widetilde{\mu}_{22}$
  & $\widetilde{\mu}_{31}$ & $\widetilde{\mu}_{32}$ & $\widetilde{\mu}_{41}$ \\
  \hline
  20  & -2.37  & -33.66&  11.14&   -4.10  & 19.07  & -0.68&   3.91\\
  40  & -10.75 & -46.42  &14.62&   -8.18  & 26.44  & -5.17   &6.28\\
  60  &-0.29 &  -58.24  &19.64  & -5.81   &32.92 &  0.97 &   6.32\\
  80  &-1.69 &  -67.11  &22.50  & -7.15 &  37.96  & 0.35&    7.45\\
  100 &-1.31  & -74.91 & 25.17  & -7.79  & 42.36  & 0.72  &  8.24\\
  110 &-1.33   &-78.51 & 26.38  & -8.15  & 44.39  & 0.77   & 8.64\\
  120 &-1.75 &  -81.49 & 27.42  & -8.54  & 46.20 &  0.65  &  9.03\\
  140 &-1.75 &  -81.49 & 27.42  & -8.54  & 46.20 &  0.65  &  9.03\\
  160 &-1.75 &  -81.49 & 27.42  & -8.54  & 46.20 &  0.65  &  9.03\\
 \hline
\end{tabular}

\medskip

\begin{tabular}{c|c|c|c|c|c|c}
  \hline
  % after \\: \hline or \cline{col1-col2} \cline{col3-col4} ...
  Trunc. & $\widetilde{\mu}_{42}$ & $\widetilde{\mu}_{51}$ & $\widetilde{\mu}_{52}$  & $BIC^{*}$ & Time & Iter. \\
  \hline
  20  & 3.71  &  -23.13 & 2.97  &  -3625.80 &   101& 2083\\
  40  & 4.30   & -31.60&  9.27  &  -4247.02  &  185& 2067\\
  60  &6.82  &  -40.17  &2.52  &  -4754.54   & 273& 2050\\
  80  &7.72    &-46.23&  3.84   & -5209.68&    322 &1816\\
  100 &8.68  &  -51.63  &3.88 &   -5634.52 &   329& 1449\\
  110 &9.10  &  -54.11&  4.05  &  -5839.09 &   440& 1776\\
  120 &9.41   & -56.30 & 4.38  &  -6019.10&    461 &1688\\
  140   &9.41   & -56.30 & 4.38  &  -6019.10&    461 &1688\\
  160 &9.41   & -56.30 & 4.38  &  -6019.10&    461 &1688\\
 \hline
\end{tabular}
\end{small}
\end{center}
\end{table}

\medskip

\begin{table}[!ht]
  \centering
  \caption{The maximum likelihood estimators for the \textit{large} group under the
  \textit{Gaussian} model}\label{tab:Gaussianlarge}
\begin{center}
\begin{small}
\renewcommand{\arraystretch}{1}
\begin{tabular}{c|c|c|c|c|c|c|c}
  \hline
  % after \\: \hline or \cline{col1-col2} \cline{col3-col4} ...
  Trunc. & $\widetilde{\mu}_{11}$  & $\widetilde{\mu}_{12}$ & $\widetilde{\mu}_{21}$ & $\widetilde{\mu}_{22}$
  & $\widetilde{\mu}_{31}$ & $\widetilde{\mu}_{32}$ & $\widetilde{\mu}_{41}$ \\
  \hline
  20  &-19.04 & -22.88 & 5.37 &   -8.26 &  15.82 &  -10.82 & 3.84\\
  40  &-29.85  &-29.93 & 6.54   & -12.38 & 20.99 &  -17.32 & 5.62\\
  60  &-15.90 & -49.41  &14.07  & -9.87  & 32.64 &  -7.20  & 5.30\\
  80  &-41.34&  -43.55 & 9.73   & -17.35 & 30.42  & -23.86  &7.93\\
  100 &-66.69 & -8.40   &-3.88  & -21.92&  9.43    &-42.24&  8.53\\
  110 &-40.91  &-57.46 & 14.17  & -18.57 & 39.32&   -22.74&  8.83\\
  120 &-32.30 & -65.98 & 17.67   &-16.67  &44.17  & -16.68  &8.38\\
  140 &-32.30 & -65.98 & 17.67   &-16.67  &44.17  & -16.68  &8.38\\
  160 &-32.30 & -65.98 & 17.67   &-16.67  &44.17  & -16.68  &8.38\\
 \hline
\end{tabular}

\medskip

\begin{tabular}{c|c|c|c|c|c|c}
  \hline
  % after \\: \hline or \cline{col1-col2} \cline{col3-col4} ...
  Trunc. & $\widetilde{\mu}_{42}$ & $\widetilde{\mu}_{51}$ & $\widetilde{\mu}_{52}$  & $BIC^{*}$ & Time & Iter. \\
  \hline
  20  &1.42 &   -18.06 & 15.31 &  -3540.51 &   155& 2075\\
  40  &1.52 &   -23.59 & 23.97 &  -4159.47 &   259& 1824\\
  60  &4.80  &  -39.19&  13.02 &  -4665.18 &   274& 1295\\
  80  &2.35 &   -34.35&  33.21   &-5118.90 &   300 &1044\\
  100 &-3.59 &  -6.20 &  53.12 &  -5542.62 &   978& 2753\\
  110 &4.04&    -45.42&  32.97   &-5746.72  &  449& 1143\\
  120 &5.65 &   -52.16  &26.15&   -5926.64&    509& 1172\\
  140 &5.65 &   -52.16  &26.15&   -5926.64&    509& 1172\\
  160 &5.65 &   -52.16  &26.15&   -5926.64&    509& 1172\\
 \hline
\end{tabular}
\end{small}
\end{center}
\end{table}

\medskip

\begin{table}[!ht]
  \centering
  \caption{The maximum likelihood estimators for the \textit{large} group under the
  \textit{Kotz $T=2$} model}\label{tab:T2large}
\begin{center}
\begin{small}
\renewcommand{\arraystretch}{1}
\begin{tabular}{c|c|c|c|c|c|c|c}
  \hline
  % after \\: \hline or \cline{col1-col2} \cline{col3-col4} ...
  Trunc. & $\widetilde{\mu}_{11}$  & $\widetilde{\mu}_{12}$ & $\widetilde{\mu}_{21}$ & $\widetilde{\mu}_{22}$
  & $\widetilde{\mu}_{31}$ & $\widetilde{\mu}_{32}$ & $\widetilde{\mu}_{41}$ \\
  \hline
  20  &-21.67 & -21.97 & 4.83 &   -9.01 &  15.40 &  -12.56 & 4.10\\
  40  & -36.15&  -24.57 & 4.23 &   -13.84  &17.94 &  -21.68 & 6.00\\
  60  &-32.77&  -42.35  &10.19  & -14.52&  29.14   &-18.44 & 6.82\\
  80  &-31.52 & -53.20 & 13.74 &  -15.19 & 36.02  & -16.98&  7.42\\
  100 & -31.91&  -61.32 & 16.28  & -16.10  &41.25   &-16.74 & 8.02\\
  110 &-38.06  &-61.70 & 15.79  & -18.09 & 41.86  & -20.66 & 8.78\\
  120 &-42.11 & -62.61 & 15.65 &  -19.44 & 42.61&   -23.16 & 9.32\\
  140 &-42.11 & -62.61 & 15.65 &  -19.44 & 42.61&   -23.16 & 9.32\\
  160 &-42.11 & -62.61 & 15.65 &  -19.44 & 42.61&   -23.16 & 9.32\\
 \hline
\end{tabular}

\medskip

\begin{tabular}{c|c|c|c|c|c|c}
  \hline
  % after \\: \hline or \cline{col1-col2} \cline{col3-col4} ...
  Trunc. & $\widetilde{\mu}_{42}$ & $\widetilde{\mu}_{51}$ & $\widetilde{\mu}_{52}$  & $BIC^{*}$ & Time & Iter. \\
  \hline
  20  & 1.13 &   -17.32 & 17.40 &  -3586.83 &   339& 2165\\
  40  &0.44 &   -19.28 & 28.94 &  -4207.80 &   476& 1594\\
  60  &2.80 &   -33.46&  26.38&   -4715.29&    519& 1133\\
  80  &4.18  &  -42.10 & 25.47  & -5170.59  &  706& 1144\\
  100 &5.12 &   -48.56 & 25.83   &-5595.75  &  1143  &  1528\\
  110 &4.74 &   -48.81 & 30.73&   -5800.52   & 1105 &   1349\\
  120 &4.58 &   -49.41 & 33.92  & -5981.01 &   1320 &   1459\\
  140 &4.58 &   -49.41 & 33.92  & -5981.01 &   1320 &   1459\\
  160 &4.58 &   -49.41 & 33.92  & -5981.01 &   1320 &   1459\\
 \hline
\end{tabular}
\end{small}
\end{center}
\end{table}

\medskip

\begin{table}[!ht]
  \centering
  \caption{The maximum likelihood estimators for the \textit{large} group under the
  \textit{Kotz $T=3$} model}\label{tab:T3large}
\begin{center}
\begin{small}
\renewcommand{\arraystretch}{1}
\begin{tabular}{c|c|c|c|c|c|c|c}
  \hline
  % after \\: \hline or \cline{col1-col2} \cline{col3-col4} ...
  Trunc. & $\widetilde{\mu}_{11}$  & $\widetilde{\mu}_{12}$ & $\widetilde{\mu}_{21}$ & $\widetilde{\mu}_{22}$
  & $\widetilde{\mu}_{31}$ & $\widetilde{\mu}_{32}$ & $\widetilde{\mu}_{41}$ \\
  \hline
  20  & -31.74 & 3.30  &  -4.16&   -9.72 &  -0.19  & -20.55 & 3.56\\
  40  &-41.25&  -18.15 & 1.70 &   -14.83  &14.14  & -25.34 & 6.17\\
  60  &-24.01 & -49.58 & 13.33&   -12.45&  33.24  & -12.39&  6.27\\
  80  &-32.14 & -54.75 & 14.17 &  -15.53 & 37.05 &  -17.29 & 7.60\\
  100 &-41.95 & -57.11 & 13.96&   -18.87 & 39.15 &  -23.43 & 8.93\\
  110 &-35.24 & -65.37 & 17.23 &  -17.55 & 44.04 &  -18.63 & 8.70\\
  120 &-39.44 & -66.43 & 17.11&   -18.97 & 44.89 &  -21.22 & 9.27\\
  140 &-39.44 & -66.43 & 17.11&   -18.97 & 44.89 &  -21.22 & 9.27\\
  160 &-39.44 & -66.43 & 17.11&   -18.97 & 44.89 &  -21.22 & 9.27\\
 \hline
\end{tabular}

\medskip

\begin{tabular}{c|c|c|c|c|c|c}
  \hline
  % after \\: \hline or \cline{col1-col2} \cline{col3-col4} ...
  Trunc. & $\widetilde{\mu}_{42}$ & $\widetilde{\mu}_{51}$ & $\widetilde{\mu}_{52}$  & $BIC^{*}$ & Time & Iter. \\
  \hline
  20  & -2.58  & 2.86  &  25.23&   -3628.05 &   73&  1494\\
  40  &-0.67  & -14.14 & 32.95 &  -4251.15 &   125 &1411\\
  60  &4.26 &   -39.27&  19.47 &  -4760.57  &  163 &1239\\
  80  &4.32 &   -43.32 & 25.97 &  -5217.61  &  197 &1110\\
  100 &3.93 &   -45.13 & 33.79 &  -5644.35 &   299 &1345\\
  110 &5.37 &   -51.75 & 28.52 &  -5849.87 &   304& 1246\\
  120 &5.21  &  -52.46 & 31.83 &  -6031.00 &   387 &1457\\
  140   &5.21  &  -52.46 & 31.83 &  -6031.00 &   387 &1457\\
    160 &5.21  &  -52.46 & 31.83 &  -6031.00 &   387 &1457\\
 \hline
\end{tabular}
\end{small}
\end{center}
\end{table}

At this point the log likelihood can be computed, then we use fminsearch for the
MLE's. The initial value for the algorithm is the sample mean of the elliptical
matrix variables $\mathbf{Y}\sim\mathcal{E}_{N-1 \times
K}(\boldsymbol{\mu}\mathbf{\Theta}^{-1/2}, \mathbf{\Sigma}\otimes \mathbf{I}_{K},h)
$. However, we need to deal with an open problem proposed by \citet{KE06},  the
relationship between the convergence and the truncation of the series.  Concretely,
how many terms we need to consider in the series (\ref{suma}) in order to reach some
fixed tolerance for convergence.  A numerical solution consists of optimising the log
likelihood, by increasing the truncation until, the MLE's and the maximum of the
function, reach an equilibrium, which depends on the standard accuracy and tolerance
of the routine fminsearch.  We tried the truncations $20, 40, 60, 80, 100, 110, 120,
140$ and $160$, and we note that after the truncation 120 the solutions stabilise.
the maximum likelihood estimators for location parameters associated with the small
and large groups under the Gaussian, Kotz $T=2$ and Kotz $T=3$ models, are summarized
in tables \ref{tab:Gaussiansmall}-\ref{tab:T3large}, respectively. Tables also show
the modified $BIC^{*}$  value, the number of iterations for obtaining the convergence
and the time in seconds for each optimisation.

\begin{figure}[h]
  \begin{center}
  \includegraphics[width=5.5cm,height=4cm]{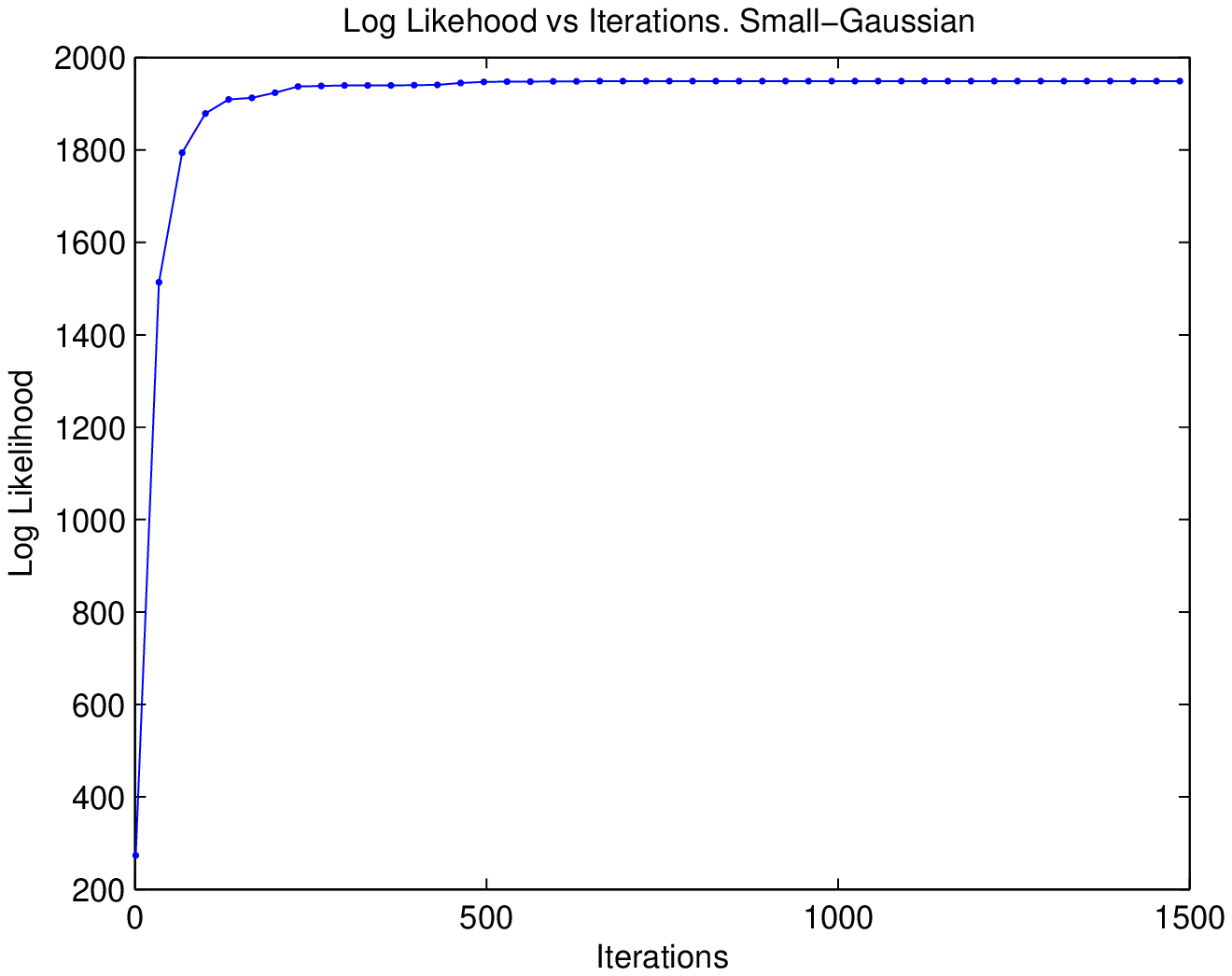}
  \includegraphics[width=5.5cm,height=4cm]{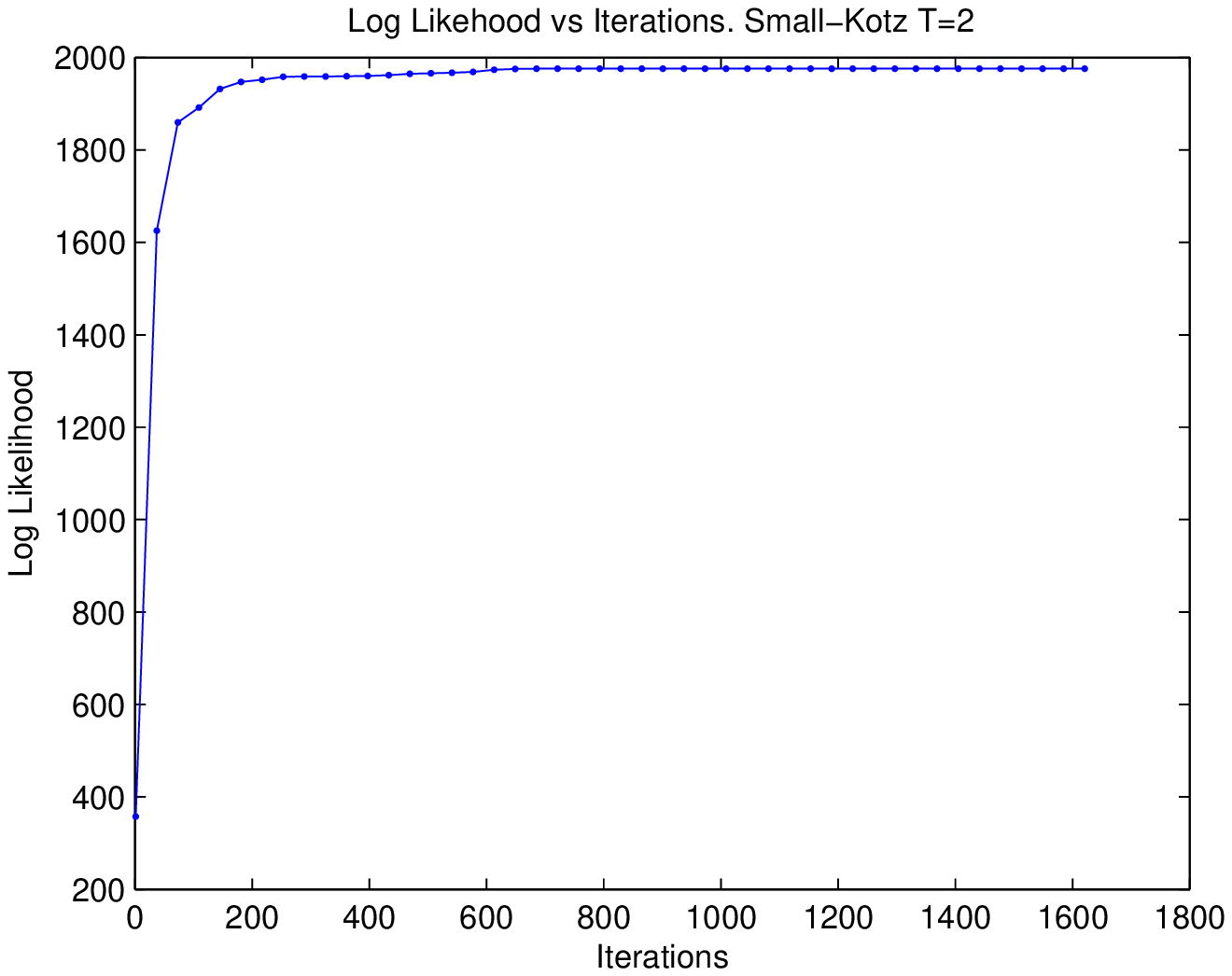}\\
  \includegraphics[width=5.5cm,height=4cm]{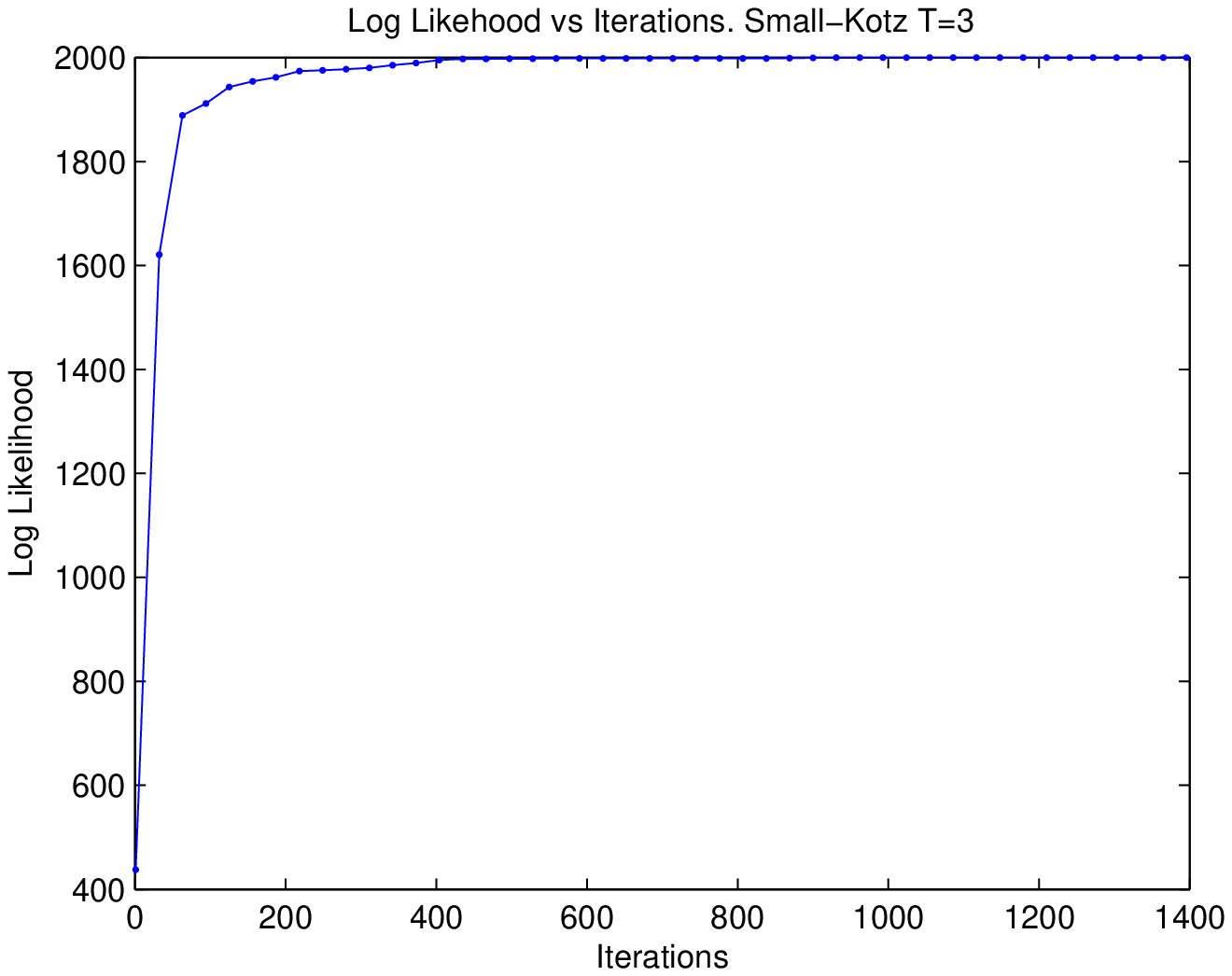}
  \includegraphics[width=5.5cm,height=4cm]{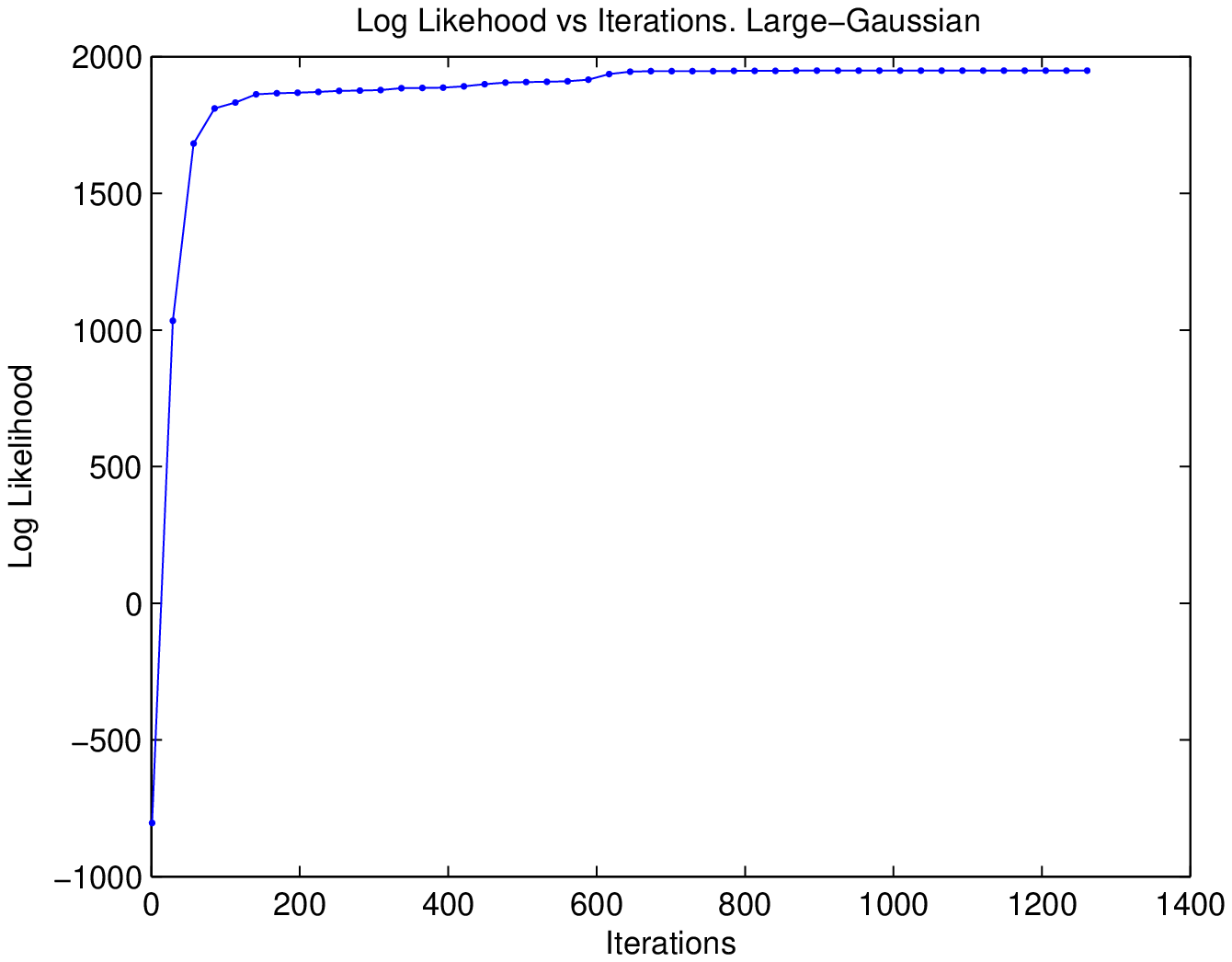}\\
  \includegraphics[width=5.5cm,height=4cm]{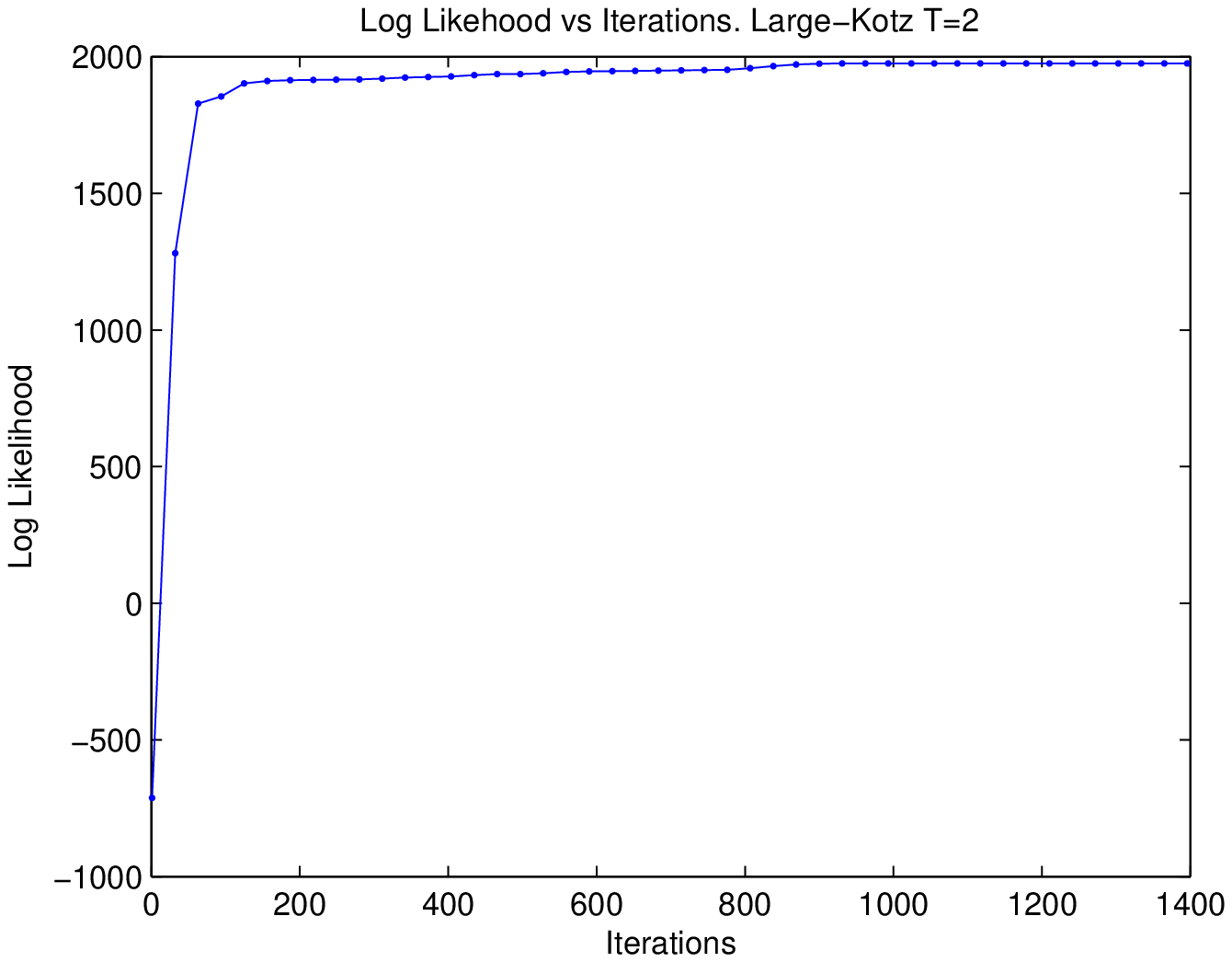}
  \includegraphics[width=5.5cm,height=4cm]{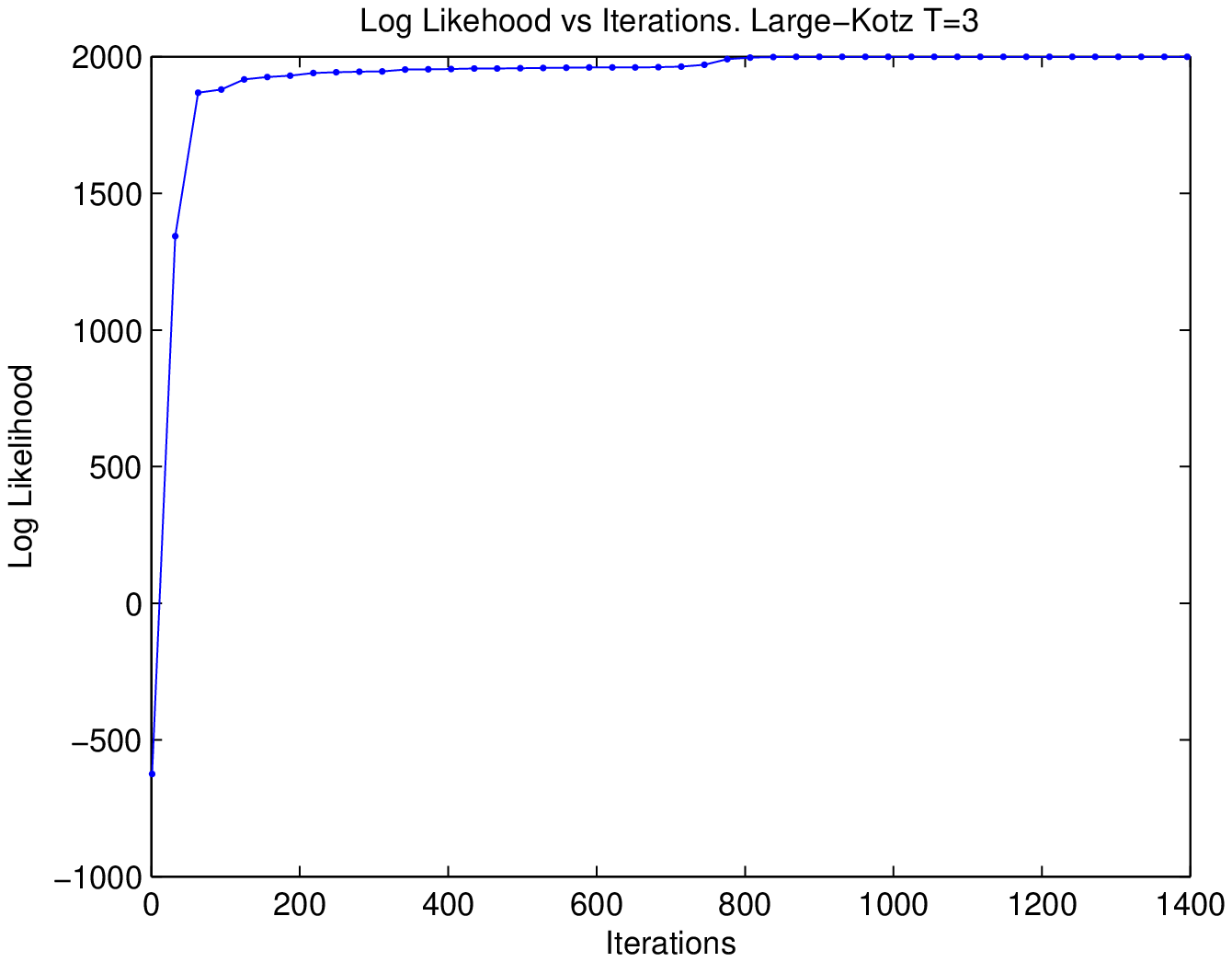}\\
  \caption{Behaviors of Log-Likelihood functions in terms of the iteration
  number of the fminsearch routine.}\label{fig:LogV}
  \end{center}
\end{figure}

The computations were performed with a processor Intel(R) Corel(TM)2 Duo CPU,
E7400@2.80GHz, and 2,96GB of RAM.

Figures \ref{fig:LogV} show the behavior of the maximum of the log likelihood when
the number of iterations is increased. In this case we use a truncation of 160, and
again, we note that the log likelihood is bounded for a very small number of
iterations in each particular model.

According to the modified BIC criterion, we can order the models in the large and
small groups as follows: (1) Kotz $T=3$, (2) Kotz $T=2$, (3) Gaussian.

\begin{figure}[h]
  \begin{center}
  \includegraphics[width=12cm,height=6cm]{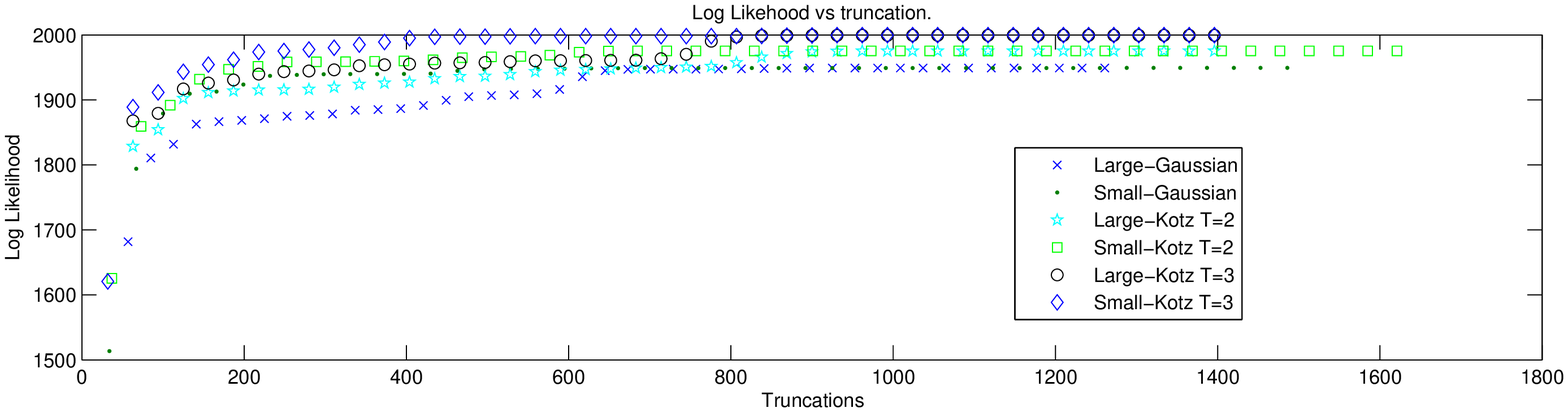}
  \caption{Comparison of Log-Likelihood in all models and groups in terms of the iteration
  number of the fminsearch routine.}\label{fig:multiple}
  \end{center}
\end{figure}

This order can be seen in figure \ref{fig:multiple}, which compares the
log-likelihood of the two groups under the three models in terms of the algorithm
iteration when the truncation is set in 160.

Modified $BIC^{*}$ of both groups shows a very strong difference (see table
\ref{table2}) between the best model (1) and the classical Gaussian (3).

In both cases, the true models of the data maybe have tails that are weighted more or
less than  Gaussian model or that the shape distribution present grater or smaller
degree of kurtosis than the Gaussian model.

\begin{rem}
We have used this example from the literature to illustrate the generalised shape
theory; moreover, based on the modified $BIC^{*}$, we found that the Kotz
distribution (with $T=3$)  is the best model in this experiment. However, suppose the
expert in the area of application knows that the landmarks have a Gaussian
distribution, then we must apply the classical theory of shape (based on normality).
Alternatively, if the expert in the application area suspects that the landmarks  do
not have a Gaussian distribution, so we can apply the generalized theory proposed
here. In this case the expert has the  necessary tools to choose an elliptical model
(as an alternative to the Gaussian distribution), according to the characteristic of
the sample which reveal and/or support a non Gaussian distribution, i.e. to select a
distribution with more or less  heavy tails, or more or less kurtosis than the
Gaussian density; among many others possible characteristics.
\end{rem}

Once the best models are selected for the small and large groups, we can test
equality in mean shape between the two independent populations. In this experiment we
have: two independent samples of 23 bones and 10 population shape parameters to
estimate for each group. Namely, if $L(\boldsymbol{\mu}_{s},\boldsymbol{\mu}_{l})$ is
the likelihood, where $\boldsymbol{\mu}_{s}$, $\boldsymbol{\mu}_{l}$,  represent the
mean shape parameters of the small and large group, respectively, then we want to
test: $H_{0}:\boldsymbol{\mu}_{s}=\boldsymbol{\mu}_{l}$ vs
$H_{a}:\boldsymbol{\mu}_{s}\neq\boldsymbol{\mu}_{l}$. Then $-2\log \Lambda=2\sup
_{H_{1}}\log L(\boldsymbol{\mu}_{s},\boldsymbol{\mu}_{l}) -2\sup _{H_{0}}\log
L(\boldsymbol{\mu}_{s}, \boldsymbol{\mu}_{l})$, and according to Wilk's theorem
$-2\log \Lambda\sim\chi^{2}_{10}$ under $H_{0}.$

Using  fminsearch with a truncation of 160  we  obtained  that:
$$
  -2\log \Lambda=2(3999.1273)-2(3990.3601)=17.5344
,
$$
this is the same result when the series were truncated at 120 and 140. Since the
p-value for the test is
$$
  P(\chi^{2}_{10}\geq 17.5344)=0.0633
$$
we have some evidence that the small and large mouse vertebrae are different in mean
shape. \citet{MD:89} studied this problem with a Gaussian model and Bookstein
coordinates (see also \citet{DM98}) and they obtained for the same test an
approximate p-value of zero ($P(\chi^{2}_{8}\geq 127.75)$. Our test also rejects the
equality of mean shape based on a better non Gaussian model but without an strong
evidence as the Gaussian model suggests.

Note that the MLE's given by tables \ref{tab:Gaussiansmall}-\ref{tab:T3small}
correspond to the matrix $\boldsymbol{\mu}$ in $\mathbf{Y}\sim\mathcal{E}_{N-1 \times
K}(\boldsymbol{\mu}, \mathbf{\Sigma}\otimes \mathbf{I}_{K},h), $ we can use this
information and the transformations

$$
   \mathbf{L}\mathbf{X}\mathbf{\Theta}^{-1/2}=\mathbf{L}\mathbf{Z}=\mathbf{Y}\Rightarrow\mathbf{V}=
    r\mathbf{W}=r\mathbf{W}(\mathbf{u}),
$$
(with $\mathbf{V}=\mathbf{YY}'$ and  $\mathbf{W} = \mathbf{V}/r$) to estimate the
different means at each step, i.e.: the original elliptical mean
$\boldsymbol{\mu}_{\mathbf{X}}$, the size-and-shape mean
$\boldsymbol{\mu}_{\mathbf{V}}$ and the shape mean $\boldsymbol{\mu}_{\mathbf{W}}$.

This example deserves a detailed study about some important facts, i.e.  the
distribution of $-2\log \Lambda$ for small samples, the truncation of the series,
global optimisation methods, etc. These problems shall be considered in a subsequent
work.

A final comment, for any elliptical model we can obtain the SVD reflection model,
however a nontrivial problem appears, the $2t$-th derivative of the generator model,
which can be seen as a partition theory problem. For the general case of a Kotz model
($s \neq 1$), and another models as Pearson II and VII, Bessel, Jensen-logistic, we
can use formulae for these derivatives given by \citet{Caro2009}. The resulting
densities have again a form of a generalised series of zonal polynomials which can be
computed efficiently after some modification of existing works for hypergeometric
series, see \citet{KE06}, thus the inference over an exact density can be performed,
avoiding the use of any asymptotic distribution, and the initial transformation
avoids the invariant polynomials of \citet{d:80}, which at present are not computable
for large degrees.

\section*{Acknowledgments}

This research work was supported  by University of Medellin (Medellin, Colombia) and
Universidad Aut\'onoma Agraria Antonio Narro (M\'{e}xico),  joint grant No. 469,
SUMMA group. Also, the first author was partially supported  by IDI-Spain, Grants No.
\ FQM2006-2271 and MTM2008-05785 and the paper was written during J. A. D\'{\i}az-
Garc\'{\i}a's stay as a visiting professor at the Department of Statistics and O. R.
of the University of Granada, Spain. Finally, F. Caro thanks to the project No. 105657
of CONACYT, M\'{e}xico.


\begin{thebibliography}{}

\bibitem[Billingsley(1986)]{bi:86}
   P. Billingsley,
   Probability and Measure,
   John Wiley \& Sons, New York, 1986.

\bibitem[Caro-Lopera \textit{et al.}(2009)]{Caro2009}
    F. J. Caro-Lopera, J. A. D\'{\i}az-Garc\'{\i}a and G. Gonz\'{a}lez-Far\'{\i}as,
    Noncentral elliptical configuration density,
    J. Multivariate Anal. 101(1) (2009), 32--43.

\bibitem[Davis(1980)]{d:80}
 A. W. Davis,
  Invariant polynomials with two matrix arguments, extending the zonal polynomials,
   in: Multivariate Analysis V, (Krishnaiah, P. R. ed.), North-Holland, 1980.

\bibitem[D\'{\i}az-Garc\'{\i}a and Guti\'errez-J\'aimez(1997)]{dg:97}
   J. A. D\'{\i}az-Garc\'{\i}a, and R. Guti\'errez-J\'aimez,
   Proof of the conjectures of H. Uhlig on the singular multivariate beta and
   the jacobian of a certain matrix transformation,
   Ann. Statist., 25, (1997) 2018-2023.

\bibitem[D\'{\i}az-Garc\'{\i}a et al.(1997)]{dgm:97}
   J. A. D\'{\i}az-Garc\'{\i}a, R. Guti\'errez-J\'aimez, and K. V. Mardia,
   Wishart and Pseudo-Wishart distributions and some applications to
   shape theory,
   J. Multivariate Anal. 63 (1997) 73-87.

\bibitem[D\'{\i}az-Garc\'{\i}a and Gonz\'{a}lez-Far\'{\i}as (2005)]{dggg05}
    J. A. D\'{\i}az-Garc\'{\i}a and G. Gonz\'{a}lez-Far\'{\i}as,
    Singular random matrix decompositions: Distributions,
    J. Multivariate Anal. 194(1) (2005), 109--122.

\bibitem[D\'{\i}az-Garc\'{\i}a and Guti\'errez-J\'aimez (2006)]{dggj06}
    J. A. D\'{\i}az-Garc\'{\i}a and R. Guti\'errez-J\'aimez,
    Wishart and Pseudo-Wishart distributions
    under elliptical laws and related distributions in the shape theory context,
    J. Stat. Plan. Inference 136(12) (2006), 4176--4193.

\bibitem[Dryden and Mardia (1998)]{DM98}
    I. L. Dryden and K.V. Mardia,
    Statistical shape analysis,
    John Wiley and Sons, Chichester, 1998.

\bibitem[Fang and Zhang (1990)]{fz:90}
    K. T. Fang,  and Y. T. Zhang, Generalized
    Multivariate Analysis,
    Science Press, Springer-Verlag, Beijing, 1990.

\bibitem[Goodall(1991)]{g:91}
    C. G. Goodall,
    Procustes methods in the statistical analysis of shape (with discussion),
    J. Roy. Statist. Soc. Ser. B, 53 (1991) 285-339.

\bibitem[Goodall and Mardia (1993)]{GM93}
    C. R. Goodall, and K. V. Mardia,
    Multivariate Aspects of Shape Theory,
    Ann.  Statist. 21 (1993) 848--866.

\bibitem[Gupta and Varga(1993)]{gv:93}
    A. K. Gupta, and T. Varga,
    Elliptically Contoured Models in Statistics,
    Kluwer Academic Publishers, Dordrecht, 1993.

\bibitem[James(1964)]{JAT64}
    A. T. James,  Distributions of matrix variate and latent roots
   derived from normal samples,
   Ann. Math. Statist. 35 (1964) 475--501.

\bibitem[Kass and Raftery (1995)]{kr:95}
    R. E. Kass, and A. E. Raftery, Bayes factor, J. Amer. Statist. Soc. 90 (1995) 773--795.

\bibitem[Khatri(1968)]{k:68}
    C. G. Khatri,
    Some results for the singular normal multivariate regression models,
    Sankhy\={a} A  30 (1968) 267-280.

\bibitem[Koev and Edelman (2006)]{KE06}
    P. Koev and A. Edelman,
    The efficient evaluation of the hypergeometric function of a matrix argument,
    Math. Comp. 75 (2006) 833--846.

\bibitem[Le and Kendall (1993)]{lk:93}
    H. L. Le, and D. G. Kendall,
    The Riemannian structure of Euclidean spaces: a novel environment for
    statistics,
    Ann.Statist. 21 (1993) 1225--1271.

\bibitem[Mardia and Dryden (1989)]{MD:89}
    K. V. Mardia and I. L. Dryden,
    The Statistical Analysis of Shape Data,
    Biometrika, 76(2) (1989) 271--281

\bibitem[Muirhead(1982)]{MR1982}
    R. J. Muirhead,
    Aspects of multivariate statistical theory,
    Wiley Series in Probability and Mathematical Statistics. John Wiley \& Sons, Inc.
    1982.

\bibitem[Raftery(1995)]{r:95}
    A. E. Raftery, Bayesian model selection in social research,
    Sociological Methodology, 25 (1995) 111--163.

\bibitem[Rao(1973)]{r:73}
    C. R. Rao,
    Linear Statistical Inference and its Applications (2nd ed.),
    John  Wiley \& Sons, New York, 1973.

\bibitem[Rissanen(1978)]{ri:78}
    J. Rissanen, Modelling by shortest data description,
    Automatica, 14 (1978) 465--471.

\bibitem[Uhlig(1994)]{u:94}
    H. Uhlig, On singular Wishart and singular multivariate Beta distributions,
    Ann. Statist. 22 (1994) 395-405.

\bibitem[Yang and Yang(2007)]{YY07}
    Ch. Ch. Yang and Ch. Ch. Yang,
    Separating latent classes by information criteria,
    J. Classification 24 (2007) 183--203.

\end{thebibliography}
\end{document}